\documentclass[a4paper,leqno,draft]{article}
\usepackage{amssymb,amsmath,amsfonts,amsthm,
mathrsfs}
\newtheorem{theorem}{Theorem}[section]
\newtheorem{corollary}[theorem]{Corollary}
\newtheorem{lemma}[theorem]{Lemma}
\newtheorem{proposition}[theorem]{Proposition}
\newtheorem{definition}[theorem]{Definition}
\theoremstyle{definition}
\newtheorem{remark}[theorem]{Remark}
\newtheorem{example}[theorem]{Example}

\newcommand{\wt}[1]{\widetilde{#1}}

\newcommand{\Cinf}{\ensuremath{\mathcal{C}^\infty}}
\newcommand{\Cinfc}{\ensuremath{\mathcal{C}^\infty_{\text{c}}}}
\newcommand{\D}{\ensuremath{{\cal D}}}
\renewcommand{\S}{\mathscr{S}}
\newcommand{\E}{\ensuremath{{\cal E}}}

\newcommand{\LL}{\mathcal{L}}

\newcommand{\mb}[1]{\ensuremath{\mathbb{#1}}}
\newcommand{\N}{\mb{N}}

\newcommand{\R}{\mb{R}}
\newcommand{\C}{\mb{C}}

\newcommand{\G}{\ensuremath{{\cal G}}}
\newcommand{\Gt}{\ensuremath{{\cal G}_\tau}}

\newcommand{\Gc}{\ensuremath{{\cal G}_\mathrm{c}}}
\newcommand{\Gcinf}{\ensuremath{{\cal G}^\infty_\mathrm{c}}}
\newcommand{\GS}{\G_{{\, }\atop{\hskip-4pt\scriptstyle\S}}\!}
\newcommand{\EM}{\ensuremath{{\cal E}_{M}}}

\newcommand{\Et}{\ensuremath{{\cal E}_{\tau}}}

\newcommand{\EMinf}{\ensuremath{{\cal E}^\infty_{M}}}

\newcommand{\Nt}{\ensuremath{{\cal N}_{\tau}}}
\newcommand{\Neg}{\mathcal{N}}

\newcommand{\Ginf}{\ensuremath{\G^\infty}}

\newcommand{\lara}[1]{\langle #1 \rangle}
\newcommand{\WF}{\mathrm{WF}}

\newcommand{\singsupp}{\mathrm{sing\, supp}}
\newcommand{\supp}{\mathrm{supp}}

\newcommand{\zs}{\setminus 0}
\newcommand{\CO}[1]{\ensuremath{T^*(#1) \zs}}
\newcommand{\ssc}{\mathrm{sc}}

\newcommand{\Ellsc}{\mathrm{Ell}_\ssc}

\newfont{\bigmath}{cmr12 at 13pt}
\newcommand{\PPsi}{\bigmath{\symbol{9}}}
\newcommand{\Oprop}[1]{{\ }_{\mathrm{pr}}^{\hphantom{m}}\text{\PPsi}_{\ssc}^{ #1}}

\newfont{\grecomath}{cmmi12 at 15pt}
\newcommand{\Bigmu}{\text{\grecomath{\symbol{22}}}}

\newcommand{\val}{\mathrm{v}} 
\newcommand{\esp}{\mathrm{e}}

\renewcommand{\d}{\ensuremath{\partial}}
\newcommand{\diff}[1]{\frac{d}{d#1}}

\newfont{\bl}{msbm10 scaled \magstep2}

\newcommand{\beq}{\begin{equation}}
\newcommand{\eeq}{\end{equation}}

\newcommand{\notmid}{\mid\kern-0.5em\not\kern0.5em}

\newcommand{\eps}{\varepsilon}

\newcommand{\Om}{\Omega}

\newcommand{\compl}[1]{{#1}^{\mathrm{c}}}

\newcommand{\Syscu}{{\wt{\underline{\mathcal{S}}}}_{\,\ssc}}

\newcommand{\mO}{\mathcal{O}}
\newcommand{\mM}{\mathcal{M}}
\newcommand{\mF}{\mathcal{F}}
\newcommand{\mP}{\mathcal{P}}
\newcommand{\MPhi}{\mathcal{M}_\Phi}

\newcommand{\rfunc}{{\rm{basic}}}

\newcommand{\dslash}{d\hspace{-0.4em}{ }^-\hspace{-0.2em}}

\begin{document}

\title{{\bf Generalized Oscillatory Integrals\\ and Fourier Integral Operators }}
\author{Claudia Garetto\footnote{Supported by FWF (Austria), grant P16820-N04 and TWF (Tyrol), grant
UNI-0404/305.}\\[0.1cm]
Institut f\"ur Grundlagen der Bauingenieurwissenschaften\\
Leopold-Franzens-Universit\"at Innsbruck\\
\texttt{claudia@mat1.uibk.ac.at}\\
\ \\
G\"{u}nther H\"{o}rmann\footnote{Supported by FWF (Austria), grant P16820-N04.}\\[0.1cm]
Fakult\"at f\"ur Mathematik\\
Universit\"at Wien\\
\texttt{guenther.hoermann@univie.ac.at}\\
\  \\
Michael Oberguggenberger\footnote{Partially supported by FWF (Austria), grant Y237}\\[0.1cm]
Institut f\"ur Grundlagen der Bauingenieurwissenschaften\\
Leopold-Franzens-Universit\"at Innsbruck\\
\texttt{michael.oberguggenberger@uibk.ac.at}
}
\maketitle
\begin{abstract}
In this article, a theory of generalized oscillatory integrals (OIs) is developed whose phase functions
as well as amplitudes may be generalized functions of Colombeau type. Based on this, generalized
Fourier integral operators (FIOs) acting on Colombeau algebras are defined. This is motivated by the
need of a general framework for partial differential operators with non-smooth coefficients and
distribution data. The mapping properties of these FIOs are studied, as is microlocal Colombeau
regularity for OIs and the influence of the FIO action on generalized wave front sets.
\end{abstract}

\setcounter{section}{-1}
\section{Introduction}

This article is part of a program that seeks to solve linear partial differential equations with non-smooth
coefficients and strongly irregular data and study the qualitative properties of the solutions.
While a well established theory with powerful analytic methods is available in the case of operators with
(relatively) smooth coefficients \cite{Hoermander:V1-4}, many models from physics involve non-smooth
variations of the physical parameters and consequently require partial differential operators where
the smoothness assumption on the coefficients is dropped. Typical examples are equations that describe the
propagation of elastic waves in discontinuous media with point sources or stationary solutions of
such equations with strongly singular potential. In such cases, the theory of distribution
does not provide a general framework in which solutions exist because of the structural restraint in
dealing with nonlinear operations (cf. \cite{HdH:01, HS:68, O:92}), as is the product of a discontinuous
function with the prospective solution.

An alternative framework is provided by the theory of Colombeau algebras of generalized functions
\cite{Colombeau:85, GKOS:01, O:92}.
In this setting, multiplication of distributions is possible and generalized solutions can be obtained that
solve the equations in a strict differential-algebraic sense. Interpreting the non-smooth coefficients and data
as elements of the Colombeau algebra, existence and uniqueness has been established
for many classes of equations by now
\cite{Biagioni:90, BO:92, BO:92b, CO:90, GH:03, HO:03, LO:91, O:88, O:92, OR:98a, OR:98b, Ste:98}.
In order to study the regularity of solutions,
microlocal techniques have to be introduced into this setting, in particular, pseudodifferential
operators with generalized amplitudes and generalized wave front sets. This has been done in the papers
\cite{GGO:03, GH:05, GH:05b, Hoermann:99, GH:04, HK:01, HOP:05, NPS:98, Pilipovic:94},
what concerns elliptic equations and hypoellipticity.

As in the classical case, Fourier integral operators arise prominently in the study of solvability
of hyperbolic equations, regularity of solutions and the inverse problem (determining the non-smooth coefficients
from the data -- an important problem in geophysics \cite{dHSto:02}). In the case of differential operators
with coefficients belonging to the Colombeau algebras, this leads to Fourier integral operators with
generalized amplitudes {\em and} generalized phase functions. The purpose of this paper is to develop the
theory of that type of Fourier integral operators and to derive first results on propagation of singularities.

We begin with the following observation. Suppose generalized Fourier integral operators have been defined as acting
on a Colombeau algebra (as will be done is this paper). Evaluating the result at a point produces a
map from the Colombeau algebra into the ring of generalized constants $\wt{\C}$, that is, an element of the
dual of the Colombeau algebra. In this way, the notion of the dual of a Colombeau algebra enters, that is,
the space of $\widetilde{{\mathbb C}}$-linear maps which are continuous with respect to the so-called sharp topology.
Thus regularity not only of Colombeau generalized functions but also of the elements of the dual space
is to be investigated. Within the Colombeau algebra ${\mathcal G}(\Omega)$ ($\Omega$ an open subset of
${\mathbb R}^n$), regularity theory is based on
the subalgebra ${\mathcal G}^\infty(\Omega)$ whose intersection with ${\mathcal D}'(\Omega)$
coincides with $\Cinf(\Omega)$. An element of the dual can be regular in more subtle ways -- it may be
defined by an element of ${\mathcal G}(\Omega)$ or by an element of ${\mathcal G}^\infty(\Omega)$.
Thus for elements of the dual, two different notions of singularity arise: the notion of ${\mathcal G}$-singular
support and the notion of ${\mathcal G}^\infty$-singular support (and similarly for the wave front sets).

Having said this, we can now describe the contents of the paper in more detail. Section 1 serves to collect material
from Colombeau theory that we need. In particular, we recall topological notions, generalized symbols
and various tools for studying regularity. Further, the ${\mathcal G}(\Omega)$- and
${\mathcal G}^\infty(\Omega)$-wave front set of a functional on the Colombeau algebra is introduced.
In Section 2, we develop the foundations for generalized Fourier integral operators -- oscillatory integrals
with generalized phase functions. As in the classical case, a generalized phase function is homogeneous
of degree one in its second variable. The classical condition that the gradient should not vanish has to be replaced by
invertibility of the norm of the gradient as a Colombeau generalized function.
Generalized oscillatory integrals are then supplemented by an additional parameter in Section 3,
leading to the notion of a Fourier integral operator with generalized amplitude and phase function.
We study the mapping properties of such operators on Colombeau algebras. As has been noticed
already in the elliptic theory \cite{GGO:03, HO:03}, two asymptotic scales are required with respect to regularity
theory using ${\mathcal G}^\infty(\Omega)$: the usual scale defining the representatives of the
elements of the Colombeau algebra and the so called slow scale. We show that Fourier integral operators
with slow scale phase function and regular amplitude map ${\mathcal G}^\infty(\Omega)$ into itself.
Section 3 also contains an example indicating how such operators arise from first order hyperbolic equations
with non-smooth coefficients.
Section 4 is devoted to investigating the functionals that are given by
generalized oscillatory integrals on the Colombeau algebra in more detail. We study the regions on
which the norm of the gradient of the phase function is not invertible and its complement.
Both regions come in two different versions, depending on what asymptotic scale is chosen (normal scale or slow scale),
which in turn correspond to ${\mathcal G}$- or ${\mathcal G}^\infty$-regularity. We find bounds on the wave front set
of these functionals, again with respect to the two notions of regularity. In the case of classical
phase functions, these bounds reduce to the classical ones involving the conic support of the
amplitude and the critical set of the phase function. In the generalized case,
this condition can only be formulated by a more complicated condition of non-invertibility.
We show how this condition of non-invertibility can be used to compute the generalized wave front
set of the kernel of the Fourier integral operator arising from first order hyperbolic equations.

\section{Basic notions: Colombeau and duality theory}
\label{section_basic}
This section gives some background of Colombeau and duality theory for the techniques used in the sequel of the current paper. As main sources we refer to \cite{Garetto:05b, Garetto:05a, GGO:03, GH:05, GKOS:01}.

\subsection{Nets of complex numbers}
Before dealing with the major points of the Colombeau construction we begin by recalling some definitions concerning elements of $\mathbb{C}^{(0,1]}$.

A net $(u_\eps)_\eps$ in $\C^{(0,1]}$ is said to be \emph{strictly nonzero} if there exist $r>0$ and $\eta\in(0,1]$ such that $|u_\eps|\ge \eps^r$ for all $\eps\in(0,\eta]$.\\
The regularity issues discussed in Sections 3 and 4 will make use of the following concept of \emph{slow scale net (s.s.n)}. A slow scale net is a net $(r_\eps)_\eps\in\C^{(0,1]}$ such that 
\[
\forall q\ge 0\, \exists c_q>0\, \forall\eps\in(0,1]\qquad\qquad\qquad\qquad |r_\eps|^q\le c_q\eps^{-1}.
\]
A net $(u_\eps)_\eps$ in $\C^{(0,1]}$ is said to be \emph{slow scale-strictly nonzero} is there exist a slow scale net $(s_\eps)_\eps$ and $\eta\in(0,1]$ such that $|u_\eps|\ge 1/s_\eps$ for all $\eps\in(0,\eta]$.

\subsection{$\wt{\C}$-modules of generalized functions based on a locally convex topological vector space $E$}
\label{subsection_G_E}
The most common algebras of generalized functions of Colombeau type as well as the spaces of generalized symbols we deal with are introduced and investigated under a topological point of view by referring to the following models.

Let $E$ be a locally convex topological vector space topologized through the family of seminorms $\{p_i\}_{i\in I}$. The elements of 
\[
\begin{split} 
\mM_E &:= \{(u_\eps)_\eps\in E^{(0,1]}:\, \forall i\in I\,\, \exists N\in\N\quad p_i(u_\eps)=O(\eps^{-N})\, \text{as}\, \eps\to 0\},\\
\mM^\ssc_E &:=\{(u_\eps)_\eps\in E^{(0,1]}:\, \forall i\in I\,\, \exists (\omega_\eps)_\eps\, \text{s.s.n.}\quad p_i(u_\eps)=O(\omega_\eps)\, \text{as}\, \eps\to 0\},\\
\mM^\infty_E &:=\{(u_\eps)_\eps\in E^{(0,1]}:\, \exists N\in\N\,\, \forall i\in I\quad p_i(u_\eps)=O(\eps^{-N})\, \text{as}\, \eps\to 0\},\\
\Neg_E &:= \{(u_\eps)_\eps\in E^{(0,1]}:\, \forall i\in I\,\, \forall q\in\N\quad p_i(u_\eps)=O(\eps^{q})\, \text{as}\, \eps\to 0\},
\end{split}
\]
 
are called $E$-moderate, $E$-moderate of slow scale type, $E$-regular and $E$-negligible, respectively. We define the space of \emph{generalized functions based on $E$} as the factor space $\G_E := \mM_E / \Neg_E$. 

The ring of \emph{complex generalized numbers}, denoted by $\wt{\C}$, is obtained by taking $E=\C$. $\wt{\C}$ is not a field since by Theorem 1.2.38 in \cite{GKOS:01} only the elements which are strictly nonzero (i.e. the elements which have a representative strictly nonzero) are invertible and vice versa. Note that all the representatives of $u\in\wt{\C}$ are strictly nonzero once we know that there exists at least one which is strictly nonzero. When $u$ has a representative which is slow scale-strictly nonzero we say that it is \emph{slow scale-invertible}.

For any locally convex topological vector space $E$ the space $\G_E$ has the structure of a $\wt{\C}$-module. The ${\C}$-module $\G^\ssc_E:=\mM^\ssc_E/\Neg_E$ of \emph{generalized functions of slow scale type} and the $\wt{\C}$-module $\Ginf_E:=\mM^\infty_E/\Neg_E$ of \emph{regular generalized functions} are subrings of $\G_E$ with more refined assumptions of moderateness at the level of representatives. We use the notation $u=[(u_\eps)_\eps]$ for the class $u$ of $(u_\eps)_\eps$ in $\G_E$. This is the usual way adopted in the paper to denote an equivalence class.

The family of seminorms $\{p_i\}_{i\in I}$ on $E$ determines a \emph{locally convex $\wt{\C}$-linear} topology on $\G_E$ (see \cite[Definition 1.6]{Garetto:05b}) by means of the \emph{valuations}
\[
\val_{p_i}([(u_\eps)_\eps]):=\val_{p_i}((u_\eps)_\eps):=\sup\{b\in\R:\qquad p_i(u_\eps)=O(\eps^b)\, \text{as $\eps\to 0$}\}
\] 
and the corresponding \emph{ultra-pseudo-seminorms} $\{\mP_i\}_{i\in I}$. For the sake of brevity we omit to report definitions and properties of valuations and ultra-pseudo-seminorms in the abstract context of $\wt{\C}$-modules. Such a theoretical presentation can be found in \cite[Subsections 1.1, 1.2]{Garetto:05b}. We recall that on $\wt{\C}$ the valuation and the ultra-pseudo-norm obtained through the absolute value in $\C$ are denoted by $\val_{\wt{\C}}$ and $|\cdot|_{\esp}$ respectively. Concerning the space $\Ginf_E$ of regular generalized functions based on $E$ the moderateness properties of $\mM_E^\infty$ allows to define the valuation 
\[
\val^\infty_E ((u_\eps)_\eps):=\sup\{b\in\R:\, \forall i\in I\qquad p_i(u_\eps)=O(\eps^b)\, \text{as $\eps\to 0$}\}
\]
which extends to $\Ginf_E$ and leads to the ultra-pseudo-norm $\mP^\infty_E(u):=\esp^{-\val_E^\infty(u)}$.

The Colombeau algebra $\G(\Om)=\EM(\Om)/\Neg(\Om)$ can be obtained as a ${\wt{\C}}$-module of $\G_E$-type by choosing $E=\E(\Om)$. Topologized through the family of seminorms $p_{K,i}(f)=\sup_{x\in K, |\alpha|\le i}|\partial^\alpha f(x)|$ where $K\Subset\Om$, the space $\E(\Om)$ induces on $\G(\Om)$ a metrizable and complete locally convex $\wt{\C}$-linear topology  which is determined by the ultra-pseudo-seminorms $\mP_{K,i}(u)=\esp^{-\val_{p_{K,i}}(u)}$. $\G(\Om)$ is continuously embedded in each $(\G_{\mathcal{C}^k(\Om)},\{\mP_{K,k}(u)\}_{K\Subset\Om})$ since $\EM(\Om)\cap\Neg_{\mathcal{C}^k(\Om)}\subseteq\EM(\Om)\cap\Neg_{\mathcal{C}^0(\Om)}=\Neg(\Om)$ and the topology on $\G(\Om)$ is finer than the topology induced by any $\G_{\mathcal{C}^k(\Om)}$ on $\G(\Om)$.
From a structural point of view $\Om\to\G(\Om)$ is a fine sheaf of differential algebras on $\R^n$.

The Colombeau algebra $\Gc(\Om)$ of generalized functions with compact support is topologized by means of a strict inductive limit procedure. More precisely, setting $\G_K(\Om):=\{u\in\Gc(\Om):\, \supp\, u\subseteq K\}$ for $K\Subset\Om$, $\Gc(\Om)$ is the strict inductive limit of the sequence of locally convex topological $\wt{\C}$-modules $(\G_{K_n}(\Om))_{n\in\N}$, where $(K_n)_{n\in\N}$ is an exhausting sequence of compact subsets of $\Om$ such that $K_n\subseteq K_{n+1}$. We recall that the space $\G_K(\Om)$ is endowed with the topology induced by $\G_{\mathcal{D}_{K'}(\Om)}$ where $K'$ is a compact subset containing $K$ in its interior. In detail we consider on $\G_K(\Om)$ the ultra-pseudo-seminorms $\mP_{\G_K(\Om),n}(u)=\esp^{-\val_{K,n}(u)}$. Note that the valuation $\val_{K,n}(u):=\val_{p_{K',n}}(u)$ is independent of the choice of $K'$ when acts on $\G_K(\Om)$. 

Regularity theory in the Colombeau context as initiated in \cite{O:92} is based on the subalgebra $\Ginf(\Om)$ of all elements $u$ of $\G(\Om)$ having a representative $(u_\eps)_\eps$ belonging to the set
\[
\EM^\infty(\Om):=\{(u_\eps)_\eps\in\E[\Om]:\ \forall K\Subset\Om\, \exists N\in\N\, \forall\alpha\in\N^n\quad \sup_{x\in K}|\partial^\alpha u_\eps(x)|=O(\eps^{-N})\, \text{as $\eps\to 0$}\}.
\]
$\Ginf(\Om)$ can be seen as the intersection $\cap_{K\Subset\Om}\Ginf(K)$, where $\Ginf(K)$ is the space of all $u\in\G(\Om)$ having a representative $(u_\eps)_\eps$ satisfying the condition: $\exists N\in\N$ $\forall\alpha\in\N^n$,\ $\sup_{x\in K}|\partial^\alpha u_\eps(x)|=O(\eps^{-N})$. The ultra-pseudo-seminorms $\mP_{\Ginf(K)}(u):=\esp^{-\val_{\Ginf(K)}}$, where 
\[
\val_{\Ginf(K)}:=\sup\{b\in\R:\, \forall\alpha\in\N^n\quad \sup_{x\in K}|\partial^\alpha u_\eps(x)|=O(\eps^b)\}
\]
equip $\Ginf(\Om)$ with the topological structure of a \emph{Fr\'echet $\wt{\C}$-module}.

Finally, let us consider the algebra $\Gcinf(\Om):=\Ginf(\Om)\cap\Gc(\Om)$. On $\Ginf_K(\Om):=\{u\in\Ginf(\Om):\, \supp\, u\subseteq K\}$ with $K\Subset\Om$, we define the ultra-pseudo-norm $\mP_{\Ginf_K(\Om)}(u)=\esp^{-\val^\infty_K(u)}$ where $\val^\infty_K(u):=\val^\infty_{\mathcal{D}_{K'}(\Om)}(u)$ and $K'$ is any compact set containing $K$ in its interior. At this point, given an exhausting sequence $(K_n)_n$ of compact subsets of $\Om$, the strict inductive limit procedure equips $\Gcinf(\Om)=\cup_n \Ginf_{K_n}(\Om)$ with a complete and separated locally convex $\wt{\C}$-linear topology.  

\subsection{Topological dual of a Colombeau algebra}
A duality theory for $\wt{\C}$-modules had been developed in \cite{Garetto:05b} in the framework of topological and locally convex topological $\wt{\C}$-modules. Starting from an investigation of $\LL(\G,\wt{\C})$, the $\wt{\C}$-module of all $\wt{\C}$-linear and continuous functionals on $\G$, it provides the theoretical tools for dealing with the topological duals of the Colombeau algebras $\Gc(\Om)$ and $\G(\Om)$. In the paper $\LL(\G(\Om),\wt{\C}$ and $\LL(\Gc(\Om),\wt{\C})$ are endowed with the \emph{topology of uniform convergence on bounded subsets}. This is determined by the ultra-pseudo-seminorms 
\[
\mP_{B^\circ}(T)=\sup_{u\in B}|T(u)|_\esp,
\]
where $B$ is varying in the family of all bounded subsets of $\G(\Om)$ and $\Gc(\Om)$ respectively. From general results concerning the relation between boundedness and ultra-pseudo-seminorms in the context of locally convex topological $\wt{\C}$-modules we have that $B\subseteq\G(\Om)$ is bounded if and only if for all $K\Subset\Om$ and $i\in\N$ there exists a constant $C>0$ such that $\mP_{K,i}(u)\le C$ for all $u\in B$. In particular the strict inductive limit structure of $\Gc(\Om)$ yields that $B\subseteq\Gc(\Om)$ is bounded if and only if it is contained in some $\G_K(\Om)$ and bounded there if and only if 
\[
\exists K\Subset\Om\, \forall n\in\N\, \exists C>0\, \forall u\in B\ \ \mP_{\G_K(\Om),n}(u)\le C.
\]
For the choice of topologies illustrated in this section Theorem 3.1 in \cite{Garetto:05a} shows the following chains of continuous embeddings:
\beq
\label{chain_1}
\Ginf(\Om)\subseteq\G(\Om)\subseteq\LL(\Gc(\Om),\wt{\C}),
\eeq
\beq
\label{chain_2}
\Gcinf(\Om)\subseteq\Gc(\Om)\subseteq\LL(\G(\Om),\wt{\C}),
\eeq  
\beq
\label{chain_3}
\LL(\G(\Om),\wt{\C})\subseteq\LL(\Gc(\Om),\wt{\C}).
\eeq  
In \eqref{chain_1} and \eqref{chain_2} the inclusion in the dual is given via integration $\big(u\to\big( v\to\int_\Om u(x)v(x)dx\big)\big)$ (for definitions and properties of the integral of a Colombeau generalized functions see \cite{GKOS:01}) while the embedding in \eqref{chain_3} is determined by the inclusion $\Gc(\Om)\subseteq\G(\Om)$. Since $\Om\to\LL(\Gc(\Om),\wt{\C})$ is a sheaf we can define the \emph{support of a functional $T$} (denoted by $\supp\, T$). In analogy with distribution theory $\LL(\G(\Om),\wt{\C})$ from Theorem 1.2 in \cite{Garetto:05a} we have that $\LL(\G(\Om),\wt{\C})$ can be identified with the set of functionals in $\LL(\Gc(\Om),\wt{\C})$ having compact support. 

By \eqref{chain_1} it is meaningful to measure the regularity of a functional in $\LL(\Gc(\Om),\wt{\C})$ with respect to the algebras $\G(\Om)$ and $\Ginf(\Om)$. We define the \emph{$\G$-singular support} of $T$ (${\rm{singsupp}}_\G\, T$) as the complement of the set of all points $x\in\Om$ such that the restriction of $T$ to some open neighborhood $V$ of $x$ belongs to $\G(V)$. Analogously replacing $\G$ with $\Ginf$ we introduce the notion of \emph{$\Ginf$-singular support} of $T$ denoted by ${\rm{singsupp}}_{\Ginf} T$. This investigation of regularity is connected with the notions of generalized wave front sets considered in Subsection \ref{sub_sec_micro} and will be focused on the functionals in $\LL(\Gc(\Om),\wt{\C})$ and $\LL(\G(\Om),\wt{\C})$ which have a ``basic'' structure. In detail, we say that $T\in\LL(\Gc(\Om),\wt{\C})$ is ${\rfunc}$ if there exists a net $(T_\eps)_\eps\in\D'(\Om)^{(0,1]}$ fulfilling the following condition: for all $K\Subset\Om$ there exist $j\in\N$, $c>0$, $N\in\N$ and $\eta\in(0,1]$ such that
\[
\forall f\in\D_K(\Om)\, \forall\eps\in(0,\eta]\qquad\quad
|T_\eps(f)|\le c\eps^{-N}\sup_{x\in K,|\alpha|\le j}|\partial^\alpha f(x)|
\]
and $Tu=[(T_\eps u_\eps)_\eps]$ for all $u\in\Gc(\Om)$.\\
In the same way a functional $T\in\LL(\G(\Om),\wt{\C})$ is said to be $\rfunc$ if there exists a net  $(T_\eps)_\eps\in\E'(\Om)^{(0,1]}$ such that there exist $K\Subset\Om$, $j\in\N$, $c>0$, $N\in\N$ and $\eta\in(0,1]$ with the property 
\[
\forall f\in\Cinf(\Om)\, \forall\eps\in(0,\eta]\qquad\quad
|T_\eps(f)|\le c\eps^{-N}\sup_{x\in K,|\alpha|\le j}|\partial^\alpha f(x)|
\]
and $Tu=[(T_\eps u_\eps)_\eps]$ for all $u\in\G(\Om)$.\\
Clearly the sets of $\rfunc$ functionals are $\wt{\C}$-linear subspaces of $\LL(\Gc(\Om),\wt{\C})$ and $\LL(\G(\Om),\wt{\C})$ respectively. In addition if $T$ is a $\rfunc$ functional in $\LL(\Gc(\Om),\wt{\C})$ and $u\in\Gc(\Om)$ then $uT\in\LL(\G(\Om),\wt{\C})$ is $\rfunc$. We recall that nets $(T_\eps)_\eps$ which define basic maps as above were already considered in \cite{Delcroix:05,DelSca:00} with slightly more general notions of moderateness and different choices of notations and language.


\subsection{Generalized symbols}
\label{subsec_gen_symb}
For the convenience of the reader we recall a few basic notions concerning the sets of symbols employed in the course of the paper. 
\paragraph{Definitions.}
Let $\Om$ be an open subset of $\R^n$, $m\in\R$ and $\rho,\delta\in[0,1]$. $S^m_{\rho,\delta}(\Om\times\R^p)$ denotes the set of symbols of order $m$ and type $(\rho,\delta)$ as introduced by H\"ormander in \cite{Hoermander:71}. The subscript $(\rho,\delta)$ is omitted when $\rho=1$ and $\delta=0$. If $V$ is an open conic set of $\Om\times\R^{p}$ we define $S^m_{\rho,\delta}(V)$ as the set of all $a\in\Cinf(V)$ such that for all $K\Subset V$, 
\[
\sup_{(x,\xi)\in K^{c}}\lara{\xi}^{-m+\rho|\alpha|-\delta|\beta|}|\partial^\alpha_\xi\partial^\beta_x a(x,\xi)|<\infty,
\]
where $K^{c}:=\{(x,t\xi):\, (x,\xi)\in K\ t\ge 1\}$.
We also make use of the space $S^1_{\rm{hg}}(\Om\times\R^p\setminus 0)$ of all $a\in S^1(\Om\times\R^p\setminus 0)$ homogeneous of degree $1$ in $\xi$. Note that the assumption of homogeneity allows to state the defining conditions above in terms of the seminorms
\[
\sup_{x\in K,\xi\in\R^p\setminus 0}|\xi|^{-1+\alpha}|\partial^\alpha_\xi\partial^\beta_x a(x,\xi)|
\]
where $K$ is any compact subset of $\Om$.

The space of generalized symbols $\wt{\mathcal{S}}^m_{\rho,\delta}(\Om\times\R^p)$ is the $\wt{\C}$-module of $\G_E$-type obtained by taking $E=S^{m}_{\rho,\delta}(\Om\times\R^p)$ equipped with the family of seminorms
\[
|a|^{(m)}_{\rho,\delta,K,j}=\sup_{x\in K,\xi\in\R^n}\sup_{|\alpha+\beta|\le j}|\partial^\alpha_\xi\partial^\beta_x a(x,\xi)|\lara{\xi}^{-m+\rho|\alpha|-\delta|\beta|},\qquad\quad K\Subset\Om,\, j\in\N.
\]
The valuation corresponding to $|\cdot|^{(m)}_{\rho,\delta,K,j}$ gives the ultra-pseudo-seminorm $\mP^{(m)}_{\rho,\delta,K,j}$. $\wt{\mathcal{S}}^m_{\rho,\delta}(\Om\times\R^p)$ topologized through the family of ultra-pseudo-seminorms $\{\mP^{(m)}_{\rho,\delta,K,j}\}_{K\Subset\Om,j\in\N}$ is a Fr\'echet $\wt{\C}$-module. In analogy with $\wt{\mathcal{S}}^m_{\rho,\delta}(\Om\times\R^p)$ we use the notation $\wt{\mathcal{S}}^m_{\rho,\delta}(V)$ for the $\wt{\C}$-module $\G_{S^m_{\rho,\delta}(V)}$.

$\wt{\mathcal{S}}^m_{\rho,\delta}(\Om_x\times\R^p_\xi)$ has the structure of a sheaf with respect to $\Om$. So it is meaningful to talk of the support with respect to $x$ of a generalized symbol $a$ ($\supp_x\, a$). In particular, when $a\in\wt{\mathcal{S}}^m_{\rho,\delta}(\Om'_x\times\Om_y\times\R^p)$ we have the notions of support with respect to $x$ (${\rm{supp}}_x a$) and support with respect to $y$ (${\rm{supp}}_y a$).\\
We define the \emph{conic support} of $a\in\wt{\mathcal{S}}^m_{\rho,\delta}(\Om\times\R^p)$ (${\rm{cone\, supp}}\, a$) as the complement of the set of points $(x_0,\xi_0)\in\Om\times\R^p$ such that there exists a relatively compact open neighborhood $U$ of $x_0$, a conic open neighborhood $\Gamma$ of $\xi_0$ and a representative $(a_\eps)_\eps$ of $a$ satisfying the condition
\beq
\label{cond_conic_supp}
\forall\alpha\in\N^p\, \forall\beta\in\N^n\, \forall q\in\N\quad \sup_{x\in U,\xi\in\Gamma}\lara{\xi}^{-m+\rho|\alpha|-\delta|\beta|}|\partial^\alpha_\xi\partial^\beta_x a_\eps(x,\xi)|=O(\eps^q)\, \text{as $\eps\to 0$}.
\eeq
By definition ${\rm{cone\, supp}}\, a$ is a closed conic subset of $\Om\times\R^p$. The generalized symbol $a$ is $0$ on $\Om\setminus\pi_x(\rm{cone\, supp}\, a)$.

\paragraph{Regular symbols.}

The space of \emph{regular symbols} $\wt{\mathcal{S}}^m_{\rho,\delta,\rm{rg}}(\Om\times\R^p)$ as introduced in \cite{GGO:03} can be topologized as a locally convex topological $\wt{\C}$-module by observing that it coincides with $\cap_{K\Subset\Om}\wt{\mathcal{S}}^m_{\rho,\delta,\rm{rg}}(K\times\R^p)$, where $\wt{\mathcal{S}}^m_{\rho,\delta,\rm{rg}}(K\times\R^p)$ is the set of all $a\in\wt{\mathcal{S}}^m_{\rho,\delta}(\Om\times\R^p)$ such that there exists a representative $(a_\eps)_\eps$ fulfilling the following property:
\beq
\label{cond_reg_N}
\exists N\in\N\ \forall j\in\N\qquad\quad |a_\eps|^{(m)}_{\rho,\delta,K,j}=O(\eps^{-N})\ \text{as $\eps\to 0$}.
\eeq
On $\wt{\mathcal{S}}^m_{\rho,\delta,\rm{rg}}(K\times\R^p)$ we define the valuation $\val^{(m)}_{\rho,\delta,K;{\rm{rg}}}$ given, at the level of representatives, by 
\[
\sup\{b\in\R:\, \forall j\in\N\qquad |a_\eps|^{(m)}_{\rho,\delta,K,j}=O(\eps^{b})\ \text{as $\eps\to 0$}\},
\]
and the corresponding ultra-pseudo-seminorm $$\mP^{(m)}_{\rho,\delta,K;{\rm{rg}}}(a)=\esp^{-\val^{(m)}_{\rho,\delta,K;{\rm{rg}}}(a)}.$$ $\wt{\mathcal{S}}^m_{\rho,\delta,\rm{rg}}(K\times\R^p)$ is endowed with the locally convex $\wt{\C}$-linear topology determined by the usual ultra-pseudo-seminorms on $\wt{\mathcal{S}}^m_{\rho,\delta}(\Om\times\R^p)$ and by $\mP^{(m)}_{\rho,\delta,K;{\rm{rg}}}$. We equip $\wt{\mathcal{S}}^m_{\rho,\delta,\rm{rg}}(\Om\times\R^p)$ with the initial topology for the injections $\wt{\mathcal{S}}^m_{\rho,\delta,\rm{rg}}(\Om\times\R^p)\to\wt{\mathcal{S}}^m_{\rho,\delta,\rm{rg}}(K\times\R^p)$. This topology is given by the family of ultra-pseudo-seminorms $\{\mP^{(m)}_{\rho,\delta,K;{\rm{rg}}}\}_{K\Subset\Om}$ and is finer than the topology induced by $\wt{\mathcal{S}}^m_{\rho,\delta}(\Om\times\R^p)$ on $\wt{\mathcal{S}}^m_{\rho,\delta,\rm{rg}}(\Om\times\R^p)$. Indeed, for all $a\in\wt{\mathcal{S}}^m_{\rho,\delta,\rm{rg}}(\Om\times\R^p)$ one has
\beq
\label{topology_reg_comp}
\mP^{(m)}_{\rho,\delta,K,j}(a)\le \mP^{(m)}_{\rho,\delta,K;{\rm{rg}}}(a).
\eeq
\paragraph{Slow scale symbols.}
In the paper the classes of the factor space $\G^{\,\ssc}_{S^m_{\rho,\delta}(\Om\times\R^p)}$ are called \emph{generalized symbols of slow scale type}. $\G^{\,\ssc}_{S^m_{\rho,\delta}(\Om\times\R^p)}$ is included in $\wt{\mathcal{S}}^m_{\rho,\delta,{\rm{rg}}}(\Om\times\R^p)$ and equipped with the topology induced by $\wt{\mathcal{S}}^m_{\rho,\delta,{\rm{rg}}}(\Om\times\R^p)$. Substituting $S^m_{\rho,\delta}(\Om\times\R^p)$ with $S^m_{\rho,\delta}(V)$ we obtain the set $\G^{\,\ssc}_{S^m_{\rho,\delta}(V)}$ of slow scale symbols on the open set $V\subseteq\Om\times(\R^p\setminus 0)$.

\paragraph{Generalized symbols of order $-\infty$.}
Different notions of regularity are related to the sets $\wt{\mathcal{S}}^{-\infty}(\Om\times\R^p)$ and $\wt{\mathcal{S}}^{-\infty}_{\rm{rg}}(\Om\times\R^p)$ of generalized symbols of order $-\infty$.\\ The space $\wt{\mathcal{S}}^{-\infty}(\Om\times\R^p)$ of generalized symbols of order $-\infty$ is defined as the $\wt{\C}$-module $\G_{S^{-\infty}(\Om\times\R^p)}$. Its elements are equivalence classes $a$ whose representatives $(a_\eps)_\eps$ have the property $|a_\eps|^{(m)}_{K,j}=O(\eps^{-N})$ as $\eps\to 0$, where $N$ depends on the order $m$ of the symbol, on the order $j$ of the derivatives and on the compact set $K\subseteq\Om$. In analogy with the construction of $\wt{\mathcal{S}}^m_{\rm{rg}}(\Om\times\R^p)$ the space $\wt{\mathcal{S}}^{-\infty}_{\rm{rg}}(\Om\times\R^p)$ of \emph{regular symbols of order $-\infty$} is introduced as $\cap_{K\Subset\Om}\wt{\mathcal{S}}^{-\infty}_{\rm{rg}}(K\times\R^p)$, where $\wt{\mathcal{S}}^{-\infty}_{\rm{rg}}(K\times\R^p)$ is the set of all $a\in\wt{\mathcal{S}}^{-\infty}(\Om\times\R^p)$ such that there exists a representative $(a_\eps)_\eps$ satisfying the condition
\[
\exists N\in\N\, \forall m\in\R\, \forall j\in\N\qquad\quad |a_\eps|^{(m)}_{K,j}=O(\eps^{-N})\ \text{as $\eps\to 0$}.
\]
\paragraph{Symbols of refined order.} In Section 4 we will make use of the sets $\wt{\mathcal{S}}^{m/-\infty}_{\rho,\delta}(\Om\times\R^p)$ and $\wt{\mathcal{S}}^{m/-\infty}_{\rho,\delta,{\rm{rg}}}(\Om\times\R^p)$ of symbols of \emph{refined order} introduced in \cite{GH:05}. These arise from a finer partitioning of the classes in $\wt{\mathcal{S}}^m_{\rho,\delta}(\Om\times\R^p)$ and $\wt{\mathcal{S}}^m_{\rho,\delta,{\rm{rg}}}(\Om\times\R^p)$ respectively, obtained through a factorization with respect to the set $\Neg^{-\infty}(\Om \times\R^p):=\Neg_{S^{-\infty}(\Om\times\R^p)}$ of negligible nets. In other words if $a$ is a (regular) generalized symbol of order $m$ then for all representatives $(a_\eps)_\eps$ of $a$ we can write
\[
\kappa((a_\eps)_\eps):=(a_\eps)_\eps+\Neg^{-\infty}(\Om\times\R^p)\subseteq(a_\eps)_\eps+\Neg^{m}_{\rho,\delta}(\Om\times\R^p)=a.
\]

\paragraph{Generalized microsupports.} The $\G$- and $\Ginf$-regularity of generalized symbols on $\Om\times\R^n$ is measured in conical neighborhoods by means of the following notions of microsupports.

Let $a\in\wt{\mathcal{S}}^l_{\rho,\delta}(\Om\times\R^n)$ and $(x_0,\xi_0)\in\CO{\Om}$. The symbol $a$ is $\G$-smoothing at $(x_0,\xi_0)$ if there exist a representative $(a_\eps)_\eps$ of $a$, a relatively compact open neighborhood $U$ of $x_0$ and a conic neighborhood $\Gamma\subseteq\R^n\setminus 0$ of $\xi_0$ such that
\begin{multline}
\label{est_micro_G}
\forall m\in\R\, \forall\alpha,\beta\in\N^n\, \exists N\in\N\, \exists c>0\, \exists\eta\in(0,1]\, \forall(x,\xi)\in U\times\Gamma\, \forall\eps\in(0,\eta]\\
|\partial^\alpha_\xi\partial^\beta_x a_\eps(x,\xi)|\le c\lara{\xi}^m\eps^{-N}.
\end{multline}
The symbol $a$ is $\Ginf$-smoothing at $(x_0,\xi_0)$ if there exist a representative $(a_\eps)_\eps$ of $a$, a relatively compact open neighborhood $U$ of $x_0$, a conic neighborhood $\Gamma\subseteq\R^n\setminus 0$ of $\xi_0$ and a natural number $N\in\N$ such that  
\begin{multline}
\label{est_micro_Ginf}
\forall m\in\R\, \forall\alpha,\beta\in\N^n\, \exists c>0\, \exists\eta\in(0,1]\, \forall(x,\xi)\in U\times\Gamma\, \forall\eps\in(0,\eta]\\
|\partial^\alpha_\xi\partial^\beta_x a_\eps(x,\xi)|\le c\lara{\xi}^m\eps^{-N}.
\end{multline}
We define the \emph{$\G$-microsupport} of $a$, denoted by $\mu\, \supp_\G(a)$, as the complement in $\CO{\Om}$ of the set of points $(x_0,\xi_0)$ where $a$ is $\G$-smoothing and the \emph{$\Ginf$-microsupport} of $a$, denoted by $\mu\, \supp_{\Ginf}(a)$, as the complement in $\CO{\Om}$ of the set of points $(x_0,\xi_0)$ where $a$ is $\Ginf$-smoothing.

When $a\in\wt{\mathcal{S}}^{\,l/-\infty}_{\rho,\delta}(\Om\times\R^n)$ we denote the complements of the sets of points $(x_0,\xi_0)\in\CO{\Om}$ where \eqref{est_micro_G} and \eqref{est_micro_Ginf} hold for some representative of $a$ by $\mu_\G(a)$ and $\mu_{\Ginf}(a)$ respectively. Note that for symbols of refined order conditions \eqref{est_micro_G} and \eqref{est_micro_Ginf} do not depend on the choice of representatives. It is clear that:
\begin{itemize}
\item[(i)] if $a\in{\wt{\mathcal{S}}}^{-\infty}(\Om\times\R^n)$ then $\mu\,\supp_\G(a)=\emptyset$;
\item[(ii)] if $a\in{\wt{\mathcal{S}}}^{-\infty}_{\rm{rg}}(\Om\times\R^n)$ then $\mu\,\supp_{\Ginf}(a)=\emptyset$;
\item[(iii)] if $a\in{\wt{\mathcal{S}}}^{\,m/-\infty}_{\rho,\delta}(\Om\times\R^n)$ and $\mu_\G(a)=\emptyset$ then $a\in\wt{\mathcal{S}}^{-\infty}(\Om\times\R^n)$;
\item[(iv)] if $a\in{\wt{\mathcal{S}}}^{\,m/-\infty}_{\rho,\delta,{\rm{rg}}}(\Om\times\R^n)$ and $\mu_{\Ginf}(a)=\emptyset$ then $a\in\wt{\mathcal{S}}^{-\infty}_{\rm{rg}}(\Om\times\R^n)$;
\item[(v)] when $a$ is a classical symbol then $\mu\,\supp(a)=\mu_\G(a)=\mu_{\Ginf}(a)$;
\item[(vi)] if $a\in\wt{\mathcal{S}}^m_{\rho,\delta}(\Om\times\R^n)$ then 
\beq
\label{ess_G}
\mu\, \supp_\G(a) = \bigcap_{(a_\eps)_\eps \in a} \mu_\G\big(\,\kappa\big((a_\eps)_\eps\big)\,\big).
\eeq
and
\beq
\label{ess_Ginf}
\mu\, \supp_{\Ginf}(a) = \bigcap_{(a_\eps)_\eps \in a} \mu_{\Ginf}\big(\,\kappa\big((a_\eps)_\eps\big)\,\big).
\eeq
\end{itemize} 

\paragraph{Continuity results.}
By simple reasoning at the level of representatives one proves that the product is a continuous $\wt{\C}$-bilinear map from $\wt{\mathcal{S}}^{m_1}_{\rho_1,\delta_1}(\Om\times\R^p)\times\wt{\mathcal{S}}^{m_2}_{\rho_2,\delta_2}(\Om\times\R^p)$ into $\wt{\mathcal{S}}^{m_1+m_2}_{\rho_3,\delta_3}(\Om\times\R^p)$ with $\rho_3=\min\{\rho_1,\rho_2\}$ and $\delta_3=\max\{\delta_1,\delta_2\}$. Furthermore the derivative-map $\partial^\alpha_\xi\partial^\beta_x:\wt{\mathcal{S}}^m_{\rho,\delta}(\Om\times\R^p)\to\wt{\mathcal{S}}^{m-\rho|\alpha|+\delta|\beta|}_{\rho,\delta}(\Om\times\R^p)$ and the map
\[
\wt{\mathcal{S}}^{m_1}_{\rho_1,\delta_1}(\Om\times\R^p)\to\wt{\mathcal{S}}^{m_2}_{\rho_2,\delta_2}(\Om\times\R^p):a\to (a_\eps)_\eps+\Neg^{m_2}_{\rho_2,\delta_2}(\Om\times\R^p), 
\]
with $m_1\le m_2$, $\rho_1\ge\rho_2$, $\delta_1\le\delta_2$, are continuous.

The product between a generalized function $u(y)$ in $\Gc(\Om)$ and a generalized symbol $a(y,\xi)$ in $\wt{\mathcal{S}}^m_{\rho,\delta}(\Om\times\R^p)$ (product defined by pointwise multiplication at the level of representatives) gives an element $a(y,\xi)u(y)$ of $\wt{\mathcal{S}}^m_{\rho,\delta}(\Om\times\R^p)$. In particular the $\wt{\C}$-bilinear map
\beq
\label{bil_product}
\Gc(\Om)\times\wt{\mathcal{S}}^m_{\rho,\delta}(\Om\times\R^p)\to\wt{\mathcal{S}}^m_{\rho,\delta}(\Om\times\R^p):(u,a)\to a(y,\xi)u(y)
\eeq
is continuous. The previous results of continuity hold between spaces of regular generalized symbols and the map in \eqref{bil_product} is continuous from $\Gcinf(\Om)\times\wt{\mathcal{S}}^m_{\rho,\delta,\rm{rg}}(\Om\times\R^p)$ into $\wt{\mathcal{S}}^m_{\rho,\delta,\rm{rg}}(\Om\times\R^p)$.
\paragraph{Integration.}
If $l<-p$ each $b\in\wt{\mathcal{S}}^l_{\rho,\delta}(\Om\times\R^p)$ can be integrated on $K\times\R^p$, $K\Subset\Om$, by setting
\[
\int_{K\times\R^p}b(y,\xi)\, dy\, d\xi :=\biggl[\biggl(\int_{K\times\R^p}b_\eps(y,\xi)\, dy\, d\xi\biggr)_\eps\biggr].
\]
Moreover if $\supp_y b\Subset\Om$ we define the integral of $b$ on $\Om\times\R^p$ as
\[
\int_{\Om\times\R^p}b(y,\xi)\, dy\, d\xi :=\int_{K\times\R^p}b(y,\xi)\, dy\, d\xi,
\]
where $K$ is any compact set containing $\supp_y b$ in its interior. Integration defines a continuous $\wt{\C}$-linear functional on this space of generalized symbols with compact support in $y$.
\begin{proposition}
\label{pro_C_k}
Let $b$ be a generalized symbol with $\supp_{y}b\Subset\Om$.
\begin{trivlist}
\item[(i)]
If $b\in\wt{\mathcal{S}}^{l}_{\rho,\delta}(\Om'\times\Om\times\R^p)$ and $l+\delta k<-p$ then
\[
\int_{\Om\times\R^p}b(x,y,\xi)\, dy\, \dslash\xi :=\biggl[\biggl(\int_{K\times\R^p}b_\eps(x,y,\xi)\, dy\,\dslash\xi\biggl)_\eps\biggl],
\]
where $K$ is any compact set of $\Om$ containing $\supp_y b$ in its interior is a well-defined element of $\G_{\mathcal{C}^k(\Om')}$.
\item[(ii)] If $b\in\wt{\mathcal{S}}^{-\infty}(\Om'\times\Om\times\R^p)$ then $\int_{\Om\times\R^p}b(x,y,\xi)\, dy\, \dslash\xi\in\G(\Om')$.
\item[(iii)] If $b\in\wt{\mathcal{S}}^{-\infty}_{\rm{rg}}(\Om'\times\Om\times\R^p)$ then $\int_{\Om\times\R^p}b(x,y,\xi)\, dy\, \dslash\xi\in\Ginf(\Om')$.
\end{trivlist}
\end{proposition}
\begin{proof}
We only give the proof of the first assertion since the second and the third are immediate.

It is clear that $(b_\eps)_\eps$ is a representative of $b$ and $l+\delta k<-p$ then $v_\eps(x):=\int_{K\times\R^p}b_\eps(x,y,\xi)\,dy\,\dslash\xi$ is a net of functions in $\mathcal{C}^k(\Om')$. More precisely, for any $\alpha\in\N^n$ with $|\alpha|\le k$ and $K'\Subset\Om'$ we have that
\beq
\label{est_int_k}
\sup_{x\in K'}|\partial^\alpha v_\eps(x)|\le \sup_{x\in K',y\in K,\xi\in\R^p}\lara{\xi}^{-l-\delta|\alpha|}|\partial^\alpha_x b_\eps(x,y,\xi)|\,\int_{\R^p}\lara{\xi}^{l+\delta|\alpha|}\, \dslash\xi.
\eeq
\eqref{est_int_k} shows that $\int_{\Om\times\R^p}b(x,y,\xi)\,dy\,\dslash\xi$ is a well-defined element of $\G_{\mathcal{C}^k(\Om')}$.
\end{proof}
\begin{remark}
\label{rem_basic}
When $l+\delta k<-p$ the inequality \eqref{est_int_k} implies 
\[
\mP_{K',k}\biggl(\int_{\Om\times\R^p}b(x,y,\xi)\, dy\, \dslash\xi\biggr)\le \mP^{(l)}_{\rho,\delta,K'\times K,k}(b)
\]
and proves that integration gives a continuous map from $\wt{\mathcal{S}}^l_{\rho,\delta}(\Om'\times\Om\times\R^p)$ to $\G_{\mathcal{C}^k(\Om')}$. In particular if $b\in\wt{\mathcal{S}}^l_{\rho,\delta,{\rm{rg}}}(\Om'\times\Om\times\R^p)$ then by \eqref{topology_reg_comp} it follows that 
\[
\mP_{K',k}\biggl(\int_{\Om\times\R^p}b(x,y,\xi)\, dy\, \dslash\xi\biggr)\le \mP^{(l)}_{\rho,\delta,K'\times K;{\rm{rg}}}(b).
\]
\end{remark}
\subsection{Microlocal analysis in the Colombeau context: generalized wave front sets in $\LL(\Gc(\Om),\wt{\C})$}
\label{sub_sec_micro}
In this subsection we recall the basic notions of microlocal analysis which involve the duals of the Colombeau algebras $\Gc(\Om)$ and $\G(\Om)$ and have been developed in \cite{Garetto:06a}. In this generalized context the role which is classically played by $\S'(\R^n)$ is given to the Colombeau algebra $\GS(\R^n):=\G_{\S(\R^n)}$. $\GS(\R^n)$ is topoloziged as in Subsection \ref{subsection_G_E} and its dual $\LL(\GS(\R^n),\wt{\C})$ is endowed with the topology of uniform convergence on bounded subsets. In the sequel $\Gt(\R^n)$ denotes the Colombeau algebra of tempered generalized functions defined as the quotient $\Et(\R^n)/\Nt(\R^n)$, where $\Et(\R^n)$ is the algebra of all  \emph{$\tau$-moderate} nets $(u_\eps)_\eps\in\Et[\R^n]:=\mO_M(\R^n)^{(0,1]}$ such that 
\[
\forall \alpha\in\N^n\, \exists N\in\N\qquad \sup_{x\in\R^n}(1+|x|)^{-N}|\partial^\alpha u_\eps(x)|=O(\eps^{-N})\qquad \text{as}\ \eps\to 0
\]
and $\Nt(\R^n)$ is the ideal of all \emph{$\tau$-negligible} nets $(u_\eps)_\eps\in\Et[\R^n]$ such that
\[
\forall \alpha\in\N^n\, \exists N\in\N\, \forall q\in\N\quad \sup_{x\in\R^n}(1+|x|)^{-N}|\partial^\alpha u_\eps(x)|=O(\eps^{q})\ \text{as}\ \eps\to 0.
\]
Theorem 3.8 in \cite{Garetto:05b} shows that we have the chain of continuous embeddings
\[
\GS(\R^n)\subseteq\Gt(\R^n)\subseteq\LL(\GS(\R^n),\wt{\C}).
\]

\paragraph{The Fourier transform on $\GS(\R^n)$, $\LL(\GS(\R^n),\wt{\C})$ and $\LL(\G(\Om),\wt{\C})$.}
The Fourier transform on $\GS(\R^n)$ is defined by the corresponding transformation at the level of representatives, as follows:
\[
\mF:\GS(\R^n)\to\GS(\R^n):u\to [(\widehat{u_\eps})_\eps].
\] 
$\mF$ is a $\wt{\C}$-linear continuous map from $\GS(\R^n)$ into itself which extends to the dual in a natural way. In detail, we define the Fourier transform of $T\in\LL(\GS(\R^n),\wt{\C})$ as the functional in $\LL(\GS(\R^n),\wt{\C})$ given by
\[
\mF(T)(u)=T(\mF u).
\]
As shown in \cite[Remark 1.5]{Garetto:06a} $\LL(\G(\Om),\wt{\C})$ is embedded in $\LL(\GS(\R^n),\wt{\C})$ by means of the map
\[
\LL(\G(\Om),\wt{\C})\to\LL(\GS(\R^n),\wt{\C}):T\to \big(u\to T(({u_\eps}_{\vert_\Om})_\eps+\Neg(\Om))\big).
\]
In particular, when $T$ is a $\rfunc$ functional in $\LL(\G(\Om),\wt{\C})$ we have from \cite[Proposition 1.6, Remark 1.7]{Garetto:06a} that the Fourier transform of $T$ is the tempered generalized function obtained as the action of $T(y)$ on $e^{-iy\xi}$, i.e., $\mF(T)=T(e^{-i\cdot\xi})=(T_\eps(e^{-i\cdot\xi}))_\eps+\Nt(\R^n)$. 

\paragraph{Generalized wave front sets of a functional in $\LL(\Gc(\Om),\wt{\C})$.} The notions of $\G$-wave front set and $\Ginf$-wave front set of a functional in $\LL(\Gc(\Om),\wt{\C})$ have been introduced in \cite{Garetto:06a} as direct analogues of the distributional wave front set in \cite{Hoermander:71}. They employ a subset of the space $\G^\ssc_{S^m(\Om\times\R^n)}$ of generalized symbols of slow scale type denoted by $\Syscu^m(\Om\times\R^n)$ and introduced in \cite[Definition 1.1]{GH:05}. We refer to \cite[Definition 1.2]{GH:05} for the definition of slow scale micro-ellipticity of $a\in\Syscu^m(\Om\times\R^n)$ and to \cite{Garetto:06a} for the action of $a(x,D)\in\Oprop{m}(\Om)$ on the dual $\LL(\Gc(\Om),\wt{\C})$. We recall that $\Oprop{m}(\Om)$ denotes the set of properly supported pseudodifferential operators with symbol in $\Syscu^m(\Om\times\R^n)$.

Let $T\in\LL(\Gc(\Om),\wt{\C})$. The $\G$-wave front set of $T$ is defined as
\[
\WF_\G T:=\bigcap_{\substack{a(x,D)\in\,\Oprop{0}(\Om)\\[0.1cm] a(x,D)T\, \in\, \G(\Om)}}\hskip-5pt \compl{\Ellsc(a)}.
\]
The $\Ginf$-wave front set of $T$ is defined as
\[
\WF_{\Ginf} T:=\bigcap_{\substack{a(x,D)\in\,\Oprop{0}(\Om)\\[0.1cm] a(x,D)T\, \in\, \Ginf(\Om)}}\hskip-5pt \compl{\Ellsc(a)}.
\]
$\WF_\G T$ and $\WF_{\Ginf}T$ are both closed conic subsets of $\CO{\Om}$. As proved in \cite{Garetto:06a}, if $T$ is a basic functional in $\LL(\Gc(\Om),\wt{\C})$ then
\[
\pi_\Om(\WF_\G T)=\singsupp_\G T
\]
and
\[
\pi_{\Om}(\WF_{\Ginf}T)=\singsupp_{\Ginf} T.
\]
\paragraph{Characterization of $\WF_\G T$ and $\WF_{\Ginf} T$ when $T$ is a basic functional.} In Section \ref{section_functional} we will employ a useful characterization of the $\G$-wave front set and the $\Ginf$-wave front set valid for functionals which are basic. It involves the sets of generalized functions $\G_{\S,0}(\Gamma)$ and $\Ginf_{\S\hskip-2pt,0}(\Gamma)$, defined on the conic subset $\Gamma$ of $\R^n\setminus 0$, as follows:
\[
\G_{\S,0}(\Gamma):=\{u\in\Gt(\R^n):\ \exists (u_\eps)_\eps\in u\ \forall l\in\R\, \exists N\in\N\quad \sup_{\xi\in\Gamma}\lara{\xi}^l|u_\eps(\xi)|=O(\eps^{-N})\, \text{as $\eps\to 0$}\},
\]
\[
\Ginf_{\S\hskip-2pt,0}(\Gamma):=\{u\in\Gt(\R^n):\ \exists (u_\eps)_\eps\in u\ \exists N\in\N\, \forall l\in\R\quad \sup_{\xi\in\Gamma}\lara{\xi}^l|u_\eps(\xi)|=O(\eps^{-N})\, \text{as $\eps\to 0$}\}.
\]
Let $T\in\LL(\Gc(\Om),\wt{\C})$. Theorem 3.13 in \cite{Garetto:06a} shows that:
\begin{itemize}
\item[(i)] $(x_0,\xi_0)\not\in\WF_\G T$ if and only if there exists a conic neighborhood $\Gamma$ of $\xi_0$ and a cut-off function $\varphi\in\Cinfc(\Om)$ with $\varphi(x_0)=1$ such that $\mF(\varphi T)\in\G_{\S,0}(\Gamma)$.
\item[(ii)] $(x_0,\xi_0)\not\in\WF_{\Ginf} T$ if and only if there exists a conic neighborhood $\Gamma$ of $\xi_0$ and a cut-off function $\varphi\in\Cinfc(\Om)$ with $\varphi(x_0)=1$ such that $\mF(\varphi T)\in\Ginf_{\S\hskip-2pt,0}(\Gamma)$.
\end{itemize}
\section{Generalized oscillatory integrals: definition}
This section is devoted to a notion of oscillatory integral where both the amplitude and the phase function are generalized objects of Colombeau type. 

In the sequel $\Om$ is an arbitrary open subset of $\R^n$.
We recall that $\phi(y,\xi)$ is a \emph{phase function} on $\Om\times\R^p$ if it is a smooth function on $\Om\times\R^p\setminus 0$, real valued, positively homogeneous of degree $1$ in $\xi$ with $\nabla_{y,\xi}\phi(y,\xi)\neq 0$ for all $y\in\Om$ and $\xi\in\R^p\setminus 0$. We denote the set of all phase functions on $\Om\times\R^p$ by $\Phi(\Om\times\R^p)$ and the set of all nets in $\Phi(\Om\times\R^p)^{(0,1]}$ by $\Phi[\Om\times\R^p]$. The notations concerning classes of symbols have been introduced in Subsection \ref{subsec_gen_symb}.
\begin{definition}
\label{def_phase_moderate}
An element of $\mathcal{M}_\Phi(\Om\times\R^p)$ is a net $(\phi_\eps)_\eps\in\Phi[\Om\times\R^p]$ satisfying the conditions:
\begin{itemize}
\item[(i)] $(\phi_\eps)_\eps\in\mathcal{M}_{S^1_{\rm{hg}}(\Om\times\R^p\setminus 0)}$,
\item[(ii)] for all $K\Subset\Om$ the net $$\biggl(\inf_{y\in K,\xi\in\R^p\setminus 0}\biggl|\nabla \phi_\eps\biggl(y,\frac{\xi}{|\xi|}\biggr)\biggr|^2\biggr)_\eps$$ is strictly nonzero. 
\end{itemize}
On $\MPhi(\Om\times\R^p)$ we introduce the equivalence relation $\sim$ as follows: $(\phi_\eps)_\eps\sim(\omega_\eps)_\eps$ if and only if $(\phi_\eps-\omega_\eps)\in\Neg_{S^1_{\rm{hg}}(\Om\times\R^p\setminus 0)}$. The elements of the factor space  $$\wt{\Phi}(\Om\times\R^p):={\mathcal{M}_\Phi(\Om\times\R^p)}/{\sim}.$$
will be called \emph{generalized phase functions}. 
\end{definition}
We shall employ the equivalence class notation $[(\phi_\eps)_\eps]$ for $\phi\in\wt{\Phi}(\Om\times\R^p)$. 

When $(\phi_\eps)_\eps$ is a net of phase functions, i.e. $(\phi_\eps)_\eps\in\Phi[\Om\times\R^p]$, Lemma 1.2.1 in \cite{Hoermander:71} shows that there exists a family of partial differential operators $(L_{\phi_\eps})_\eps$ such that ${\ }^tL_{\phi_\eps}e^{i\phi_\eps}=e^{i\phi_\eps}$ for all $\eps\in(0,1]$. $L_{\phi_\eps}$ is of the form 
\beq
\label{def_L_phi_cl}
\sum_{j=1}^p a_{j,\eps}(y,\xi)\frac{\partial}{\partial_{\xi_j}} +\sum_{k=1}^n b_{k,\eps}(y,\xi)\frac{\partial}{\partial_{y_k}}+ c_\eps(y,\xi),
\eeq
where the coefficients $(a_{j,\eps})_\eps$ belong to $S^0[\Om\times\R^p]$ and $(b_{k,\eps})_\eps$, $(c_\eps)_\eps$ are elements of $S^{-1}[\Om\times\R^p]$. 
\begin{proposition}
\label{prop_operator}
If $(\phi_\eps)_\eps\in\MPhi(\Om\times\R^p)$ then $(a_{j,\eps})_\eps\in\mM_{S^0(\Om\times\R^p)}$ for all $j=1,...,p$, $(b_{k,\eps})_\eps\in\mM_{S^{-1}(\Om\times\R^p)}$ for all $k=1,...,n$, and $(c_\eps)_\eps\in\mM_{S^{-1}(\Om\times\R^p)}$.
\end{proposition} 
The proof of this proposition requires the following lemma.
\begin{lemma}
\label{lemma_1}
Let $$\varphi_{\phi_\eps}(y,\xi):=|\nabla\phi_\eps(y,\xi/|\xi|)|^{-2}.$$ If $(\phi_\eps)_\eps\in\MPhi(\Om\times\R^p)$ then $(\varphi_{\phi_\eps})_\eps\in\mM_{S^0_{\rm{hg}}(\Om\times\R^p\setminus 0)}$.
\end{lemma} 
\begin{proof}
One easily sees that $(\varphi_{\phi_\eps})_\eps$ is a net of symbols of order $0$ on $\Om\times\R^p\setminus 0$ homogeneous in $\xi$. The moderateness is obtained by combining the fact that $(\phi_\eps)_\eps\in\mM_{S^1_{\rm{hg}}(\Om\times\R^p\setminus 0)}$ with the fact that the gradient of $(\phi_\eps)_\eps$ is strictly nonzero by Definition \ref{def_phase_moderate}$(ii)$.
\end{proof}
\begin{proof}[Proof of Proposition \ref{prop_operator}]
Let $\chi\in\Cinfc(\R^p)$ such that $\chi(\xi)=1$ for $|\xi|<1/4$ and $\chi(\xi)=0$ for $|\xi|>1/2$. From the proof of Lemma 1.2.1 in \cite{Hoermander:71} we have that
\[
\begin{split}
{a}_{j,\eps}(y,\xi) &=i\,(1-\chi(\xi))\varphi_{\phi_\eps}(y,\xi){\partial_{\xi_j}\phi_\eps}(y,\xi),\\
{b}_{k,\eps}(y,\xi) &=i\,(1-\chi(\xi))|\xi|^{-2}\varphi_{\phi_\eps}(y,\xi){\partial_{y_k}\phi_\eps}(y,\xi),\\
c_\eps &=\chi(\xi)+\sum_{j=1}^p\partial_{\xi_j} {a}_{j,\eps}+\sum_{k=1}^n\partial_{y_k}{b}_{k,\eps}.
\end{split}
\]
By Lemma \ref{lemma_1} and the properties of $\chi$ it follows that $((1-\chi)\varphi_{\phi_\eps})_\eps\in\mM_{S^0(\Om\times\R^p)}$ and $((1-\chi)|\xi|^{-2}\varphi_{\phi_\eps})_\eps\in\mM_{S^{-2}(\Om\times\R^p)}$. Moreover, $(\phi_\eps)_\eps\in\mathcal{M}_{S^1_{\rm{hg}}(\Om\times\R^p\setminus 0)}$ implies that the nets  $(\partial_{\xi_j}\phi_\eps)_\eps$ and $(\partial_{y_k}\phi_\eps)_\eps$ belong to $\mathcal{M}_{S^0_{\rm{hg}}(\Om\times\R^p\setminus 0)}$ and $\mathcal{M}_{S^1_{\rm{hg}}(\Om\times\R^p\setminus 0)}$ respectively. This allows to conclude that $(a_{j,\eps})_\eps\in\mM_{S^0(\Om\times\R^p)}$, $(b_{k,\eps})_\eps\in\mM_{S^{-1}(\Om\times\R^p)}$ and $(c_\eps)_\eps\in\mM_{S^{-1}(\Om\times\R^p)}$.
\end{proof}
We proceed by comparing the families of partial differential operators $L_{\phi_\eps}$ and $L_{\omega_\eps}$ when $(\phi_\eps)_\eps\sim(\omega_\eps)_\eps$. This makes use of the following technical lemma.
\begin{lemma}
\label{lemma_1_prop}
If $(\phi_\eps)_\eps, (\omega_\eps)_\eps\in\MPhi(\Om\times\R^p)$ and $(\phi_\eps)_\eps\sim(\omega_\eps)_\eps$ then
\beq
\label{lemma_for_1}
\big(({\partial_{\xi_j}\phi_\eps})\varphi_{\phi_\eps}-({\partial_{\xi_j}\omega_\eps})\varphi_{\omega_\eps}\big)_\eps\in\Neg_{S^0_{\rm{hg}}(\Om\times\R^p\setminus 0)}
\eeq
for all $j=1,...,p$ and
\beq
\label{lemma_for_2}
\big(({\partial_{y_k}\phi_\eps}){|\xi|^{-2}\varphi_{\phi_\eps}}-({\partial_{y_k}\omega_\eps}){|\xi|^{-2}\varphi_{\omega_\eps}}\big)_\eps\in\Neg_{S^{-1}_{\rm{hg}}(\Om\times\R^p\setminus 0)}
\eeq
for all $k=1,...,n$.
\end{lemma}
\begin{proof}
The nets in \eqref{lemma_for_1} and \eqref{lemma_for_2} are of the form $a/|b|^2-c/|d|^2$ where $a,c$ are nets of moderate type, $b,d$ are $p+n$-vectors with components of moderate type and $|b|^2, |d|^2$ are strictly nonzero nets. We can write $a/|b|^2-c/|d|^2$ as
\beq
\label{frazione}
[a(d-b)\cdot(d+b)+|b|^2(a-c)]/|b|^2 |d|^2
\eeq
Since $d-b$ is a vector with negligible components, $a-c$ is a negligible net and all the other terms in \eqref{frazione} are moderate we have that $a/|b|^2-c/|d|^2$ is negligible itself.

Concerning the net in \eqref{lemma_for_1} we have that: $a=(\partial_{\xi_j}\phi_\eps)_\eps\in\mM_{S^0_{\rm{hg}}(\Om\times\R^p\setminus 0)}$, $b=(\nabla\phi_\eps(y,\xi/|\xi|))_\eps\in(\mM_{S^0_{\rm{hg}}(\Om\times\R^p\setminus 0)})^{n+p}$, $c=(\partial_{\xi_j}\omega_\eps)_\eps\in\mM_{S^0_{\rm{hg}}(\Om\times\R^p\setminus 0)}$, $d=(\nabla\omega_\eps(y,\xi/|\xi|))_\eps\in(\mM_{S^0_{\rm{hg}}(\Om\times\R^p\setminus 0)})^{n+p}$ and $d-b\in(\Neg_{S^0_{\rm{hg}}(\Om\times\R^p\setminus 0)})^{n+p}$, $a-c\in\Neg_{S^0_{\rm{hg}}(\Om\times\R^p\setminus 0)}$, ${1}/{|b|^2},{1}/{|d|^2}\in\mM_{S^0_{\rm{hg}}(\Om\times\R^p\setminus 0)}$.
Therefore, we obtain that   $((\partial_{\xi_j}\phi_\eps)\varphi_{\phi_\eps}-(\partial_{\xi_j}\omega_\eps)\varphi_{\omega_\eps})_\eps\in\Neg_{S^0_{\rm{hg}}(\Om\times\R^p\setminus 0)}$. Assertion \eqref{lemma_for_2} is proved in the same way arguing with nets of symbols of order $-1$.
\end{proof}
An inspection of the proof of Proposition \ref{prop_operator} combined with Lemma \ref{lemma_1_prop} leads to the following result.
\begin{proposition}
\label{prop_comparison_osc}
If $(\phi_\eps)_\eps, (\omega_\eps)_\eps\in\MPhi(\Om\times\R^p)$ and $(\phi_\eps)_\eps\sim(\omega_\eps)_\eps$ then
\beq
\label{L_phi_eps-L_omega_eps}
L_{\phi_\eps}-L_{\omega_\eps}=\sum_{j=1}^p a'_{j,\eps}(y,\xi)\frac{\partial}{\partial_{\xi_j}} +\sum_{k=1}^n b'_{k,\eps}(y,\xi)\frac{\partial}{\partial_{y_k}}+ c'_\eps(y,\xi),
\eeq
where $(a'_{j,\eps})_\eps\in\Neg_{S^{0}(\Om\times\R^p)}$, $(b'_{k,\eps})_\eps\in\Neg_{S^{-1}(\Om\times\R^p)}$ and $(c'_\eps)_\eps\in\Neg_{S^{-1}(\Om\times\R^p)}$ for all $j=1,...,p$ and $k=1,...,n$.
\end{proposition}
As a consequence of Propositions \ref{prop_operator} and \ref{prop_comparison_osc} we claim that any generalized phase function $\phi\in\wt{\Phi}(\Om\times\R^p)$ defines a generalized partial differential operator
\[
L_\phi(y,\xi,\partial_y,\partial_\xi)=\sum_{j=1}^p a_{j}(y,\xi)\frac{\partial}{\partial_{\xi_j}} +\sum_{k=1}^n b_{k}(y,\xi)\frac{\partial}{\partial_{y_k}}+ c(y,\xi)
\]
whose coefficients $\{a_j\}_{j=1}^p$ and $\{b_k\}_{k=1}^n$, $c$ are generalized symbols in $\wt{\mathcal{S}}^{0}(\Om\times\R^p)$ and $\wt{\mathcal{S}}^{-1}(\Om\times\R^p)$, respectively. By construction, $L_\phi$ maps $\wt{\mathcal{S}}^m_{\rho,\delta}(\Om\times\R^p)$ continuously into $\wt{\mathcal{S}}^{m-s}_{\rho,\delta}(\Om\times\R^p)$, where $s=\min\{\rho,1-\delta\}$. Hence $L^k_\phi$ is continuous from $\wt{\mathcal{S}}^m_{\rho,\delta}(\Om\times\R^p)$ to $\wt{\mathcal{S}}^{m-ks}_{\rho,\delta}(\Om\times\R^p)$.
\begin{proposition}
\label{prop_exp}
Let $\phi\in\wt{\Phi}(\Om\times\R^p)$. The exponential $$e^{i\phi(y,\xi)}$$ is a well-defined element of $\wt{\mathcal{S}}^1_{0,1}(\Om\times\R^p\setminus 0)$.
\end{proposition}
\begin{proof}
We leave it to the reader to check that if $(\phi_\eps)_\eps\in\mM_\Phi(\Om\times\R^p)$ then $(e^{i\phi_\eps(y,\xi)})_\eps\in\mM_{S^0_{0,1}(\Om\times\R^p\setminus 0)}$. When $(\phi_\eps)_\eps\sim(\omega_\eps)_\eps$, the equality
\[
e^{i\omega_\eps(y,\xi)}-e^{i\phi_\eps(y,\xi)}=e^{i\omega_\eps(y,\xi)}\big(1-e^{i(\phi_\eps-\omega_\eps)(y,\xi)}\big)=e^{i\omega_\eps(y,\xi)}\sum_{j=1}^p  e^{i(\phi_\eps-\omega_\eps)(y,\theta\xi)}\partial_{\xi_j}(\phi_\eps-\omega_\eps)(y,\theta\xi)i\xi_j,
\]
with $\theta\in(0,1)$, implies that
\beq
\label{neg_esp}
\sup_{y\in K,\xi\in\R^p\setminus 0}|\xi|^{-1}\big|e^{i\omega_\eps(y,\xi)}-e^{i\phi_\eps(y,\xi)}\big|=O(\eps^q)
\eeq
for all $q\in\N$. At this point writing $\partial^\alpha_\xi\partial^\beta_y e^{i\phi_\eps(y,\xi)}$ as
\[
\partial^\alpha_\xi\partial^\beta_y e^{i\phi_\eps(y,\xi)}\big(1-e^{i(\phi_\eps-\omega_\eps)(y,\xi)}\big)+\sum_{\alpha'<\alpha,\beta'<\beta}\binom{\alpha}{\alpha'}\binom{\beta}{\beta'}\partial^{\alpha'}_\xi\partial^{\beta'}_y e^{i\omega_\eps(y,\xi)}\big(-\partial^{\alpha-\alpha'}_\xi\partial^{\beta-\beta'}_y e^{i(\phi_\eps-\omega_\eps)(y,\xi)}\big)
\]
we obtain the characterizing estimate of a net in $\Neg_{S^1_{0,1}(\Om\times\R^p\setminus 0)}$, using \eqref{neg_esp} and the moderateness of $(e^{i\phi_\eps(y,\xi)})_\eps$.
\end{proof}
By construction of the operator $L_\phi$ the equality ${\ }^t L_\phi e^{i\phi}=e^{i\phi}$ holds in $\wt{\mathcal{S}}^1_{0,1}(\Om\times\R^p\setminus 0)$. In addition, Proposition \ref{prop_exp} and the properties of $L^k_\phi$ allow to conclude that $$e^{i\phi(y,\xi)}L^k_\phi(a(y,\xi)u(y))$$ is a generalized symbol in $\wt{\mathcal{S}}^{m-ks+1}_{0,1}(\Om\times\R^p)$ which is integrable on $\Om\times\R^p$ in the sense of Section \ref{section_basic} when $m-ks+1<-p$. From now on we assume that $\rho>0$ and $\delta<1$.
\begin{definition}
\label{def_gen_osc}
Let $\phi\in\wt{\Phi}(\Om\times\R^p)$, $a\in\wt{\mathcal{S}}^m_{\rho,\delta}(\Om\times\R^p)$ and $u\in\Gc(\Om)$. The \emph{generalized oscillatory integral}
\[
\int_{\Om\times\R^p}e^{i\phi(y,\xi)}a(y,\xi)u(y)\, dy\,\dslash\xi
\]
is defined as
\[
\int_{\Om\times\R^p}e^{i\phi(y,\xi)}L^k_\phi(a(y,\xi)u(y))\, dy\,\dslash\xi
\]
where $k$ is chosen such that $m-ks+1<-p$.  
\end{definition}
The functional 
\[
I_\phi(a):\Gc(\Om)\to\wt{\C}:u\to\int_{\Om\times\R^p}e^{i\phi(y,\xi)}a(y,\xi)u(y)\, dy\, \dslash\xi
\]
belongs to the dual $\LL(\Gc(\Om),\wt{\C})$. Indeed, by \eqref{bil_product}, the continuity of $L^k_\phi$ and of the product between generalized symbols we have that the map
\[
\Gc(\Om)\to\wt{\mathcal{S}}^{m-ks+1}_{0,1}(\Om\times\R^p):u\to e^{i\phi(y,\xi)}L^k_\phi(a(y,\xi)u(y))
\]
is continuous and thus, by an application of the integral on $\Om\times\R^p$, the resulting functional $I_\phi(a)$ is continuous.

\section{Generalized Fourier integral operators}
We now study oscillatory integrals where an additional parameter $x$, varying in an open subset $\Om'$ of $\R^{n'}$, appears in the phase function $\phi$ and in the symbol $a$. The dependence on $x$ is investigated in the Colombeau context. We denote by $\Phi[\Om';\Om\times\R^p]$ the set of all nets $(\phi_\eps)_{\eps\in(0,1]}$ of continuous functions on $\Om'\times\Om\times\R^p$ which are smooth on $\Om'\times\Om\times\R^p\setminus\{0\}$ and such that $(\phi_\eps(x,\cdot,\cdot))_\eps\in\Phi[\Om\times\R^p]$ for all $x\in\Om'$.
\begin{definition}
\label{def_phase_x_moderate}
An element of $\mM_{\Phi}(\Om';\Om\times\R^p)$ is a net $(\phi_\eps)_\eps\in\Phi[\Om';\Om\times\R^p]$ satisfying the conditions:
\begin{itemize}
\item[(i)] $(\phi_\eps)_\eps\in\mM_{S^1_{\rm{hg}}(\Om'\times\Om\times\R^p\setminus 0)}$,
\item[(ii)] for all $K'\Subset\Om'$ and $K\Subset\Om$ the net
\beq
\label{net_FIO}
\biggl(\inf_{x\in K',y\in K,\xi\in\R^p\setminus 0}\biggl|\nabla_{y,\xi} \phi_\eps\biggl(x,y,\frac{\xi}{|\xi|}\biggr)\biggr|^2\biggr)_\eps
\eeq
is strictly nonzero.
\end{itemize}
On $\mM_{\Phi}(\Om';\Om\times\R^p)$ we introduce the equivalence relation $\sim$ as follows: $(\phi_\eps)_\eps\sim(\omega_\eps)_\eps$ if and only if $(\phi_\eps-\omega_\eps)_\eps\in\Neg_{S^1_{\rm{hg}}(\Om'\times\Om\times\R^p\setminus 0)}$. The elements of the factor space 
\[
\wt{\Phi}(\Om';\Om\times\R^p):=\mM_{\Phi}(\Om';\Om\times\R^p) / \sim.
\]
are called \emph{generalized phase functions with respect to the variables in $\Om\times\R^p$}.  
\end{definition}
Proposition \ref{prop_operator} can be adapted to nets in $\mM_{\Phi}(\Om';\Om\times\R^p)$. More precisely, the operator
\beq
\label{smilla}
L_{\phi_\eps}(x;y,\xi,\partial_y,\partial_\xi)=\sum_{j=1}^p a_{j,\eps}(x,y,\xi)\frac{\partial}{\partial_{\xi_j}} +\sum_{k=1}^n b_{k,\eps}(x,y,\xi)\frac{\partial}{\partial_{y_k}}+ c_\eps(x,y,\xi)
\eeq
defined for any value of $x$ by \eqref{def_L_phi_cl}, has the property ${\ }^tL_{\phi_\eps(x,\cdot,\cdot)}e^{i\phi_\eps(x,\cdot,\cdot)}=e^{i\phi_\eps(x,\cdot,\cdot)}$ for all $x\in\Om'$ and $\eps\in(0,1]$ and its coefficients depend smoothly on $x\in\Om'$.
\begin{proposition}
\label{prop_operator_x}
If $(\phi_\eps)_\eps\in\mM_{\Phi}(\Om';\Om\times\R^p)$ then the coefficients occurring in \eqref{smilla} satisfy the following: $(a_{j,\eps})_\eps\in\mM_{S^0(\Om'\times\Om\times\R^p)}$ for all $j=1,...,p$, $(b_{k,\eps})_\eps\in\mM_{S^{-1}(\Om'\times\Om\times\R^p)}$ for all $k=1,...,n$, and $(c_\eps)_\eps\in\mM_{S^{-1}(\Om'\times\Om\times\R^p)}$.
\end{proposition}
The proof of Proposition \ref{prop_operator_x} employs the following lemma concerning basic properties of the term  $|\nabla_{y,\xi}\phi_\eps(x,y,\xi/|\xi|)|^{-2}$.  
\begin{lemma}
\label{lemma_1_x}
Let 
\beq
\label{def_varphi}
\varphi_{\phi_\eps}(x,y,\xi):=|\nabla_{y,\xi}\phi_\eps(x,y,\xi/|\xi|)|^{-2}.
\eeq
If $(\phi_\eps)_\eps\in\mM_{\Phi}(\Om';\Om\times\R^p)$ then $(\varphi_{\phi_\eps})_\eps\in\mM_{S^0_{\rm{hg}}(\Om'\times\Om\times\R^p\setminus 0)}$.
\end{lemma} 

We leave it to the reader to check that Lemma \ref{lemma_1_prop} can be stated for nets of phase functions in $(y,\xi)$ and leads to negligible nets of amplitudes in $S^0_{\rm{hg}}(\Om'\times\Om\times\R^p\setminus 0)$ and $S^{-1}_{\rm{hg}}(\Om'\times\Om\times\R^p\setminus 0)$. As a consequence we have a result on the dependence of $L_{\phi_\eps}$ on the phase function.
\begin{proposition}
\label{prop_comparison_osc_x}
If $(\phi_\eps)_\eps, (\omega_\eps)_\eps\in\mM_{\Phi}(\Om';\Om\times\R^p)$ and $(\phi_\eps)_\eps\sim(\omega_\eps)_\eps$ then 
\beq
\label{L_phi_eps-L_omega_eps_x}
L_{\phi_\eps}-L_{\omega_\eps}=\sum_{j=1}^p a'_{j,\eps}(x,y,\xi)\frac{\partial}{\partial_{\xi_j}} +\sum_{k=1}^n b'_{k,\eps}(x,y,\xi)\frac{\partial}{\partial_{y_k}}+ c'_\eps(x,y,\xi),
\eeq
where $(a'_{j,\eps})_\eps\in\Neg_{S^{0}(\Om'\times\Om\times\R^p)}$, $(b'_{k,\eps})_\eps\in\Neg_{S^{-1}(\Om'\times\Om\times\R^p)}$ and $(c'_\eps)_\eps\in\Neg_{S^{-1}(\Om'\times\Om\times\R^p)}$ for all $j=1,...,p$ and $k=1,...,n$.
\end{proposition}
Combining Propositions \ref{prop_operator_x} and \ref{prop_comparison_osc_x} yields that any generalized phase function $\phi$ in $\wt{\Phi}(\Om';\Om\times\R^p)$ defines a partial differential operator
\beq
\label{def_L_phi_amp}
L_{\phi}(x;y,\xi,\partial_y,\partial_\xi)=\sum_{j=1}^p a_{j}(x,y,\xi)\frac{\partial}{\partial_{\xi_j}} +\sum_{k=1}^n b_{k}(x,y,\xi)\frac{\partial}{\partial_{y_k}}+ c(x,y,\xi)
\eeq
with coefficients $a_j\in\wt{\mathcal{S}}^0(\Om'\times\Om\times\R^p)$, $b_k,c\in\wt{\mathcal{S}}^{-1}(\Om'\times\Om\times\R^p)$ such that ${\ }^tL_\phi e^{i\phi}=e^{i\phi}$ holds in $\wt{\mathcal{S}}^1_{0,1}(\Om'\times\Om\times\R^p\setminus 0)$.
Arguing as in Proposition \ref{prop_exp} we obtain that $e^{i\phi(x,y,\xi)}$ is a well-defined element of $\wt{\mathcal{S}}^1_{0,1}(\Om'\times\Om\times\R^p\setminus 0)$. The usual composition argument implies that the map 
\[
\Gc(\Om)\to\wt{\mathcal{S}}^{m-ks+1}_{0,1}(\Om'\times\Om\times\R^p): u\to e^{i\phi(x,y,\xi)}L^k_{\phi}(a(x,y,\xi)u(y))
\]
is continuous.
 
The oscillatory integral  
\[
I_\phi(a)(u)(x)=\int_{\Om\times\R^p}e^{i\phi(x,y,\xi)}a(x,y,\xi)u(y)\, dy\, \dslash\xi:=\int_{\Om\times\R^p}e^{i\phi(x,y,\xi)}L^k_{\phi}(a(x,y,\xi)u(y))\, dy\, \dslash\xi,
\]
where $\phi\in\wt{\Phi}(\Om';\Om\times\R^p)$ and $a\in\wt{\mathcal{S}}^m_{\rho,\delta}(\Om'\times\Om\times\R^p)$ is an element of $\wt{\C}$ for fixed $x\in\Om'$. In particular, $I_\phi(a)(u)$ is the integral on $\Om\times\R^p$ of a generalized amplitude in $\wt{\mathcal{S}}^{l}_{0,1}(\Om'\times\Om\times\R^p)$ having compact support in $y$. The order $l$ can be chosen arbitrarily low. 
\begin{theorem}
\label{theorem_map}
Let $\phi\in\wt{\Phi}(\Om';\Om\times\R^p)$, $a\in\wt{\mathcal{S}}^m_{\rho,\delta}(\Om'\times\Om\times\R^p)$ and $u\in\Gc(\Om)$. The generalized oscillatory integral
\beq
\label{oscillatory_x}
I_{\phi}(a)(u)(x)=\int_{\Om\times\R^p}e^{i\phi(x,y,\xi)}a(x,y,\xi)u(y)\, dy\,\dslash\xi  
\eeq
defines a generalized function in $\G(\Om')$ and the map
\beq
\label{def_A_fourier}
A:\Gc(\Om)\to\G(\Om'):u\to I_{\phi}(a)(u)
\eeq
is continuous.
\end{theorem}
\begin{proof}
By Proposition \ref{pro_C_k} it follows that $I_\phi(a)(u)$ is a generalized function in $\G_{\mathcal{C}^0(\Om)}$ and that for all $k\in\N$  the net
\beq
\label{repr_h}
\biggl(\int_{\Om\times\R^p}e^{i\phi_\eps(x,y,\xi)}L^h_{\phi_\eps}(a_\eps(x,y,\xi)u_\eps(y))\, dy\, \dslash\xi\biggr)_\eps
\eeq
where $sh>m+k+p+1$, belongs to $\mM_{\mathcal{C}^k(\Om)}$ and it is a representative of $I_\phi(a)(u)$. By classical arguments valid for fixed $\eps$ we know that the net given by the oscillatory integral
\[
\int_{\Om\times\R^p}e^{i\phi_\eps(x,y,\xi)}a_\eps(x,y,\xi)u_\eps(y)\, dy\,\dslash\xi
\]
is an element of $\E[\Om]$ which coincides with \eqref{repr_h} for every $k\in\N$. This means that $I_\phi(a)(u)$ is a generalized function in $\G_{\mathcal{C}^0(\Om)}$ which has a representative in $\EM(\Om)$, i.e., $I_\phi(a)(u)\in\G(\Om)$. For any $k\in\N$ the generalized function $I_\phi(a)(u)$ belongs to $\G_{\mathcal{C}^k(\Om')}$ and by Remark \ref{rem_basic} the map 
\[
\wt{\mathcal{S}}^{m-hs+1}_{0,1}(\Om'\times\Om\times\R^p)\to\G_{\mathcal{C}^k(\Om')}:e^{i\phi(x,y,\xi)}L^h_\phi(a(x,y,\xi)u(y))\to\int_{\Om\times\R^p}\hskip-10pt e^{i\phi_\eps(x,y,\xi)}L^h_{\phi_\eps}(a_\eps(x,y,\xi)u_\eps(y))\, dy\, \dslash\xi
\]
is continuous. This combined with the continuity of the map 
\[
\Gc(\Om)\to\wt{\mathcal{S}}^{m-hs+1}_{0,1}(\Om'\times\Om\times\R^p):u\to e^{i\phi(x,y,\xi)}L^h_\phi(a(x,y,\xi)u(y))  
\]
for arbitrary big $h$ proves that the map $A:u\to I_\phi(a)(u)$ is continuous from $\Gc(\Om)$ to $\G(\Om')$.
\end{proof}
The operator $A$ defined in \eqref{def_A_fourier} is called \emph{generalized Fourier integral operator} with amplitude $a\in\wt{\mathcal{S}}^m_{\rho,\delta}(\Om'\times\Om\times\R^p)$ and phase function $\phi\in\wt{\Phi}(\Om';\Om\times\R^p)$.
\begin{remark}
\label{remark_cont_a}
By continuity of the $\wt{\C}$-bilinear map $\Gc(\Om)\times\wt{\mathcal{S}}^m_{\rho,\delta}(\Om'\times\Om\times\R^p)\to\wt{\mathcal{S}}^m_{\rho,\delta}(\Om'\times\Om\times\R^p):(u,a)\to a(x,y,\xi)u(y)$ it is also clear that for fixed $u\in\Gc(\Om)$ and $\phi\in\wt{\Phi}(\Om';\Om\times\R^p)$ the map 
\[
\wt{\mathcal{S}}^{m}_{\rho,\delta}(\Om'\times\Om\times\R^p)\to\G(\Om'):a\to I_\phi(a)(u)
\]
is continuous.
\end{remark}

\begin{example}\label{G_example1}
Our outline of a basic theory of Fourier integral operators with Colombeau generalized
amplitudes and phase functions is motivated to a large extent by potential applications
in regularity theory for generalized solutions to hyperbolic partial (or pseudo-)
differential equations with distributional or Colombeau-type coefficients (or symbols)
and data (cf.\ \cite{HdH:01,LO:91,O:89}). To illustrate the typical situation we
consider here the following simple model: let $u\in\G(\R^2)$ be the solution of the
generalized Cauchy-problem
\begin{align}
\d_t u + c\,\d_x u  + b\,u &= 0\\
u \mid_{t=0} &= g,
\end{align}
where $g$ belongs to $\Gc(\R)$ and the coefficients $b$, $c\in\G(\R^2)$. 
Furthermore, $b$, $c$, as well as $\d_x c$ are
assumed to be of local $L^\infty$-log-type (concerning growth with respect to the
regularization parameter, cf.\ \cite{O:89}), $c$ being generalized real-valued and globally bounded in
addition. Let $\gamma \in \G(\R^3)$ be the unique (global) solution of the
corresponding generalized characteristic ordinary differential equation
\begin{align*}
\diff{s} \gamma(x,t;s) &= c(\gamma(x,t;s),s)\\
\gamma(x,t;t) &= x.
\end{align*}
Then $u$ is given in terms of $\gamma$ by $ u(x,t) = g(\gamma(x,t;0)) \exp(-\int_0^t
b(\gamma(x,t;r),r)\, dr)$. Writing $g$ as the inverse of its Fourier transform we
obtain the Fourier integral representation
\begin{equation}\label{hypsolu}
  u(x,t) = \iint e^{i(\gamma(x,t;0)-y) \xi}\; a(x,t,y,\xi)\, g(y)\, dy\, \dslash\xi,
\end{equation}
where $a(x,t,y,\xi) := \exp(-\int_0^t b(\gamma(x,t;r),r)\, dr)$ is a generalized
amplitude of order $0$. The phase function $\phi(x,t,y,\xi) :=
(\gamma(x,t;0)-y) \xi$ has (full) gradient
$(\d_x\gamma(x,t;0),\d_t\gamma(x,t;0),-\xi,\gamma(x,t;0)-y)$ and thus defines a generalized phase
function $\phi$. Therefore (\ref{hypsolu}) reads $u = A g$ where  $A : \Gc(\R) \to
\G(\R^2)$ is a generalized Fourier integral operator.
\end{example}

We now investigate the regularity properties of the \emph{generalized Fourier integral operator} $A$. We will prove that for appropriate generalized phase functions and generalized amplitudes, $A$ maps $\Gcinf(\Om)$ into $\Ginf(\Om')$. As in \cite{GGO:03} we consider regular amplitudes, i.e., elements $a$ of the factor space $\widetilde{\mathcal{S}}^m_{\rho,\delta,\rm{rg}}(\Om'\times\Om\times\R^p)$ whose representatives $(a_\eps)_\eps$ satisfy the condition given in \eqref{cond_reg_N} on each compact set of $\Om'\times\Om$.
However, for the phase functions the same kind of regularity assumption with respect to the parameter $\eps$ does not entail the desired mapping property.
\begin{example}
\label{example_not_reg}
Let $n=n'=p=1$ and $\Om=\Om'=\R$ and $\phi_\eps(x,y,\xi)=(x-\eps y)\xi$. Then $(\phi_\eps)_\eps\in\mM_{\Phi}(\R;\R\times\R)$ and in particular we have $N=0$ in all moderateness estimates (see Definition \ref{def_phase_x_moderate}$(i)$)) and $|\nabla_{y,\xi}\phi_\eps(x,y,\xi/|\xi|)|^2\ge\eps^2$. Choose the amplitude $a$ identically equal to $1$. The corresponding generalized operator $A$ does not map $\Gcinf(\R)$ into $\Ginf(\R)$. Indeed, for $0\neq f\in\Cinfc(\R)$ we have that
\[
A[(f)_\eps]=\biggl[\biggl(\int_{\R\times\R} e^{i(x-\eps y)\xi}f(y)\, dy\, \dslash\xi\biggr)_\eps\biggr]=[(\eps^{-1}f(x/\eps))_\eps]\in\G(\R)\setminus\Ginf(\R).
\]
\end{example}
Example \ref{example_not_reg} suggests that a stronger notion of regularity on generalized phase functions has to be designed. Such is provided by the concept of \emph{slow scale net} (s.s.n.). 
\begin{definition}
\label{def_slow_phase}
We say that $\phi\in\wt{\Phi}(\Om';\Om\times\R^p)$ is a \emph{slow scale generalized phase function in the variables of $\Om\times\R^p$} if it has a representative $(\phi_\eps)_\eps$ fulfilling the conditions
\begin{itemize}
\item[(i)] $(\phi_\eps)_\eps\in\mM^\ssc_{S^1_{\rm{hg}}(\Om'\times\Om\times\R^p\setminus 0)}$, 
\item[(ii)] for all $K'\Subset\Om'$ and $K\Subset\Om$ the net \eqref{net_FIO} is slow scale-strictly nonzero.
\end{itemize}
\end{definition}
In the sequel the set of all $(\phi_\eps)_\eps\in\Phi[\Om';\Om\times\R^p]$ fulfilling $(i)$ and $(ii)$ in Definition \ref{def_slow_phase} will be denoted by $\mM_{\Phi,\ssc}(\Om';\Om\times\R^p)$ while we use $\wt{\Phi}_\ssc(\Om';\Om\times\R^p)$ for the set of slow scale generalized functions as above. Similarly, using $\nabla_{x,y,\xi}$ in place of $\nabla_{y,\xi}$ in $(ii)$ we define the space $\wt{\Phi}_{\ssc}(\Om'\times\Om\times\R^p)$ of slow scale generalized phase functions on $\Om'\times\Om\times\R^p$.  

In the case of slow scale generalized phase functions and regular or slow scale generalized amplitudes, a careful inspection of the proofs of Proposition \ref{prop_operator_x} and Lemma \ref{lemma_1_x} leads to the following properties concerning the partial differetial operator $L_\phi$ and the Fourier integral operator $A$. We begin by observing that if $(\phi_\eps)_\eps\in\mM_{\Phi,\ssc}(\Om';\Om\times\R^p)$ then the net $(\varphi_{\phi_\eps})_\eps$ given by \eqref{def_varphi} is an element of $\mM^\ssc_{S^0_{\rm{hg}}(\Om'\times\Om\times\R^p\setminus 0)}$. Hence when $\phi\in\wt{\Phi}_{\ssc}(\Om';\Om\times\R^p)$ the operator $L_\phi$ in \eqref{def_L_phi_amp} has coefficients $a_j\in\G^\ssc_{S^0(\Om'\times\Om\times\R^p)}$ and $b_k,c\in\G^\ssc_{S^{-1}(\Om'\times\Om\times\R^p)}$. It follows that $L^h_\phi$ is continuous from $\wt{\mathcal{S}}^m_{\rho,\delta,{\rm{rg}}}(\Om'\times\Om\times\R^p)$ to $\wt{\mathcal{S}}^{m-hs}_{\rho,\delta,{\rm{rg}}}(\Om'\times\Om\times\R^p)$ and since the coefficients of $L_\phi$ are of slow scale type the inequality 
\[
\mP^{(m-hs)}_{\rho,\delta,K'\times K;{\rm{rg}}}(L^h_\phi a)\le \esp\,
\mP^{(m)}_{\rho,\delta,K'\times K;{\rm{rg}}}(a)
\]
holds for all $a$.

We will state the theorem on the regularity properties of $$A:u\to\int_{\Om\times\R^p}e^{i\phi(x,y,\xi)}a(x,y,\xi)u(y)\, dy\, \dslash\xi$$ when $\phi\in\wt{\Phi}_{\ssc}(\Om';\Om\times\R^p)$ and $a\in\wt{\mathcal{S}}^m_{\rho,\delta,\rm{rg}}(\Om'\times\Om\times\R^p)$ below. But first we observe that the product between $a(x,y,\xi)$ and $u(y)$ is a continuous $\wt{\C}$-bilinear form from $\wt{\mathcal{S}}^m_{\rho,\delta,\rm{rg}}(\Om'\times\Om\times\R^p)\times\Gcinf(\Om)$ to $\wt{\mathcal{S}}^m_{\rho,\delta,\rm{rg}}(\Om'\times\Om\times\R^p)$ and that $e^{i\phi(x,y,\xi)}\in\G^\ssc_{{S}^{1}_{0,1}(\Om'\times\Om\times\R^p\setminus 0)}$ when $\phi\in\wt{\Phi}_{\ssc}(\Om';\Om\times\R^p)$. Consequently, if $\phi\in\wt{\Phi}_{\ssc}(\Om';\Om\times\R^p)$, $a\in\wt{\mathcal{S}}^m_{\rho,\delta,\rm{rg}}(\Om'\times\Om\times\R^p)$ and $u\in\Gcinf(\Om)$ then
\[
e^{i\phi(x,y,\xi)}L^h_{\phi}(a(x,y,\xi)u(y))\in\wt{\mathcal{S}}^{m-hs+1}_{0,1,\rm{rg}}(\Om'\times\Om\times\R^p).
\]
In this situation the oscillatory integral 
\[
\int_{\Om\times\R^p}e^{i\phi(x,y,\xi)}a(x,y,\xi)u(y)\, dy\, \dslash\xi = \int_{\Om\times\R^p}e^{i\phi(x,y,\xi)}L^h_{\phi}(a(x,y,\xi)u(y))\, dy\, \dslash\xi
\]
is the integral of a generalized symbol in $\wt{\mathcal{S}}^{m-hs+1}_{0,1,\rm{rg}}(\Om'\times\Om\times\R^p)$ with compact support in $y$. 
\begin{theorem}
\label{theorem_Ginf_map}
Let $\phi\in\wt{\Phi}_{\ssc}(\Om';\Om\times\R^p)$.
\begin{itemize}
\item[(i)] If $a\in\wt{\mathcal{S}}^m_{\rho,\delta,\rm{rg}}(\Om'\times\Om\times\R^p)$ the corresponding generalized Fourier integral operator
\[
A:u\to\int_{\Om\times\R^p}e^{i\phi(x,y,\xi)}a(x,y,\xi)u(y)\, dy\, \dslash\xi
\]
maps $\Gcinf(\Om)$ continuously into $\Ginf(\Om')$.
\item[(ii)] If $a\in\wt{\mathcal{S}}^{-\infty}_{\rm{rg}}(\Om'\times\Om\times\R^p)$ then $A$ maps $\Gc(\Om)$ continuously into $\Ginf(\Om')$.
\end{itemize}
\end{theorem}
\begin{proof}
$(i)$ By Theorem \ref{theorem_map} we already know that $Au\in\G(\Om')$. A direct inspection of the representative given in \eqref{repr_h} shows that when $\phi\in\wt{\Phi}_{\ssc}(\Om';\Om\times\R^p)$ and $a\in\wt{\mathcal{S}}^m_{\rho,\delta,{\rm{rg}}}(\Om'\times\Om\times\R^p)$ the generalized function $Au$ has a representative in $\EMinf(\Om')$, i.e., $Au\in\Ginf(\Om')$. Concerning the continuity of the map $A$ we recall that $Au\in\G_{\mathcal{C}^k(\Om')}$ for each $k\in\N$ and that by Remark \ref{rem_basic}
\[
\mP_{K',k}(Au)\le \mP^{(m-hs+1)}_{0,1,K'\times K;{\rm{rg}}}\big(e^{i\phi(x,y,\xi)}L^h_\phi(a(x,y,\xi)u(y))\big),
\]
when $m-hs+1+k<-p$. Hence, the continuity of $L^h_\phi$ and of the map  $\Gcinf(\Om)\to\wt{\mathcal{S}}^m_{\rho,\delta,\rm{rg}}(\Om'\times\Om\times\R^p):u\to a(x,y,\xi)u(y)$ yields 
\begin{multline*}
\mP^{(m-hs+1)}_{0,1,K'\times K;{\rm{rg}}}\big(e^{i\phi(x,y,\xi)}L^h_\phi(a(x,y,\xi)u(y))\big)\\
\le \mP^{(1)}_{0,1,K'\times K;{\rm{rg}}}\big(e^{i\phi(x,y,\xi)}\big)\,\mP^{(m-hs)}_{\rho,\delta,K'\times K;{\rm{rg}}}(L^h_\phi(a(x,y,\xi)u(y)))\\
\le \esp^2\,\mP^{(m)}_{\rho,\delta,K'\times K;{\rm{rg}}}(a)\,\mP_{\Ginf_K(\Om)}(u)
\end{multline*}
As a consequence, since $\mP_{\Ginf(K')}(Au)=\sup_{k\in\N}\mP_{K',k}(Au)$, we may conclude that there exists a constant $C>0$ such that $\mP_{\Ginf(K')}(Au)\le C \mP_{\Ginf_K(\Om)}(u)$ for all $u$ in $\Gcinf(\Om)$ with compact support contained in $K\Subset\Om$. This proves that $A$ is continuous from $\Gcinf(\Om)$ to $\Ginf(\Om')$.

$(ii)$ Let us assume that $a\in\wt{\mathcal{S}}^{-\infty}_{\rm{rg}}(\Om'\times\Om\times\R^p)$. By Theorem \ref{theorem_map} we have that $Au\in\G(\Om)$. Since the order of $a$ is $-\infty$ the integral we deal with is absolutely convergent. Again by Remark \ref{rem_basic}, when $u\in\G_K(\Om)$ and $m+k<-p$ we obtain that 
\[
\mP_{K',k}(Au)\le \mP^{(1)}_{0,1,K'\times K;{\rm{rg}}}(e^{i\phi(x,y,\xi)})\,\mP^{(m)}_{\rho,\delta,K'\times K;{\rm{rg}}}(a)\,\mP_{\G_K(\Om),0}(u).
\]
By definition of a regular symbol of order $-\infty$ and the corresponding ultra-pseudo-seminorms, there exists a constant $C>0$ which depends only on $K'\times K$ such that $\mP^{(m)}_{\rho,\delta,K'\times K;{\rm{rg}}}(a)\le C$ for all $m$. Therefore, $Au\in\Ginf(\Om')$ and the continuity of the map $A$ is proved by 
\[
\mP_{\Ginf(K')}(Au)=\sup_{k\in\N}\mP_{K',k}(Au)\le C'\mP_{\G_K(\Om),0}(u).
\]
\end{proof}

\section{The functional $I_{\phi}(a)\in\LL(\Gc(\Om),\wt{\C})$}
\label{section_functional}
In this section we investigate the properties of the functional 
\[
I_{\phi}(a):\Gc(\Om)\to\wt{\C}:u\to\int_{\Om\times\R^p}e^{i\phi(y,\xi)}a(y,\xi)u(y)\, dy\, \dslash\xi 
\]
in more depth. Before defining specific regions depending on the generalized phase function $\phi$, we observe that any $\phi\in\wt{\Phi}(\Om\times\R^p)$ can be regarded as an element of $\G_{S^1_{\rm{hg}}(\Om\times\R^p\setminus 0)}$ and consequently $|\nabla_\xi\phi|^2\in\G_{S^0_{\rm{hg}}(\Om\times\R^p\setminus 0)}$.

Let $\Om_1$ be an open subset of $\Om$ and $\Gamma\subseteq\R^p\setminus 0$. We say that $b\in\wt{\mathcal{S}}^0(\Om\times\R^p\setminus 0)$ is \emph{invertible on $\Om_1\times\Gamma$} if for all relatively compact subsets $U$ of $\Om_1$ there exists a representative $(b_\eps)_\eps$ of $b$, a constant $r\in\R$ and $\eta\in(0,1]$ such that 
\beq
\label{est_inv_sym}
\inf_{y\in U,\xi\in\Gamma}|b_\eps(y,\xi)|\ge \eps^r
\eeq
for all $\eps\in(0,\eta]$. In an analogous way we say that $b\in\wt{\mathcal{S}}^0(\Om\times\R^p\setminus 0)$ is \emph{slow scale-invertible} on $\Om_1\times\Gamma$ if \eqref{est_inv_sym} holds with the inverse of some slow scale net $(s_\eps)_\eps$ in place of $\eps^r$. This kind of bounds from below hold for all representatives of the symbol $b$ once they are known to hold for one.

In the sequel $\pi_\Om$ denotes the projection of $\Om\times\R^p$ on $\Om$.
\begin{definition}
\label{def_C_phi}
Let $\phi\in\wt{\Phi}(\Om\times\R^p)$. We define $C_\phi\subseteq\Om\times\R^p\setminus 0$ as the complement of the set of all $(x_0,\xi_0)\in\Om\times\R^p\setminus 0$ with the property that there exist a relatively compact open neighborhood $U(x_0)$ of $x_0$ and a conic open neighborhood $\Gamma(\xi_0)\subseteq\R^p\setminus 0$ of $\xi_0$ such that $|\nabla_\xi\phi|^2$ is invertible on $U(x_0)\times\Gamma(\xi_0)$.
We set $\pi_\Om(C_{\phi})=S_{\phi}$ and $R_{\phi}=(S_{\phi})^{{\rm{c}}}$. 
\end{definition}
By construction $C_{\phi}$ is a closed conic subset of $\Om\times\R^p\setminus 0$ and $R_{\phi}\subseteq\Om$ is open. It is routine to check that the region $C_\phi$ coincides with the classical one when $\phi$ is classical.
\begin{proposition}
\label{prop_R_phi}
The generalized symbol $|\nabla_\xi\phi|^2$ is invertible on $R_\phi\times\R^p\setminus 0$. 
\end{proposition}
\begin{proof}
Let us fix a representative $(\phi_\eps)_\eps$ of $\wt{\Phi}(\Om\times\R^p)$. If $K\Subset R_{\phi}$ then $K\times\{\xi:\, |\xi|=1\}\subseteq (C_{\phi})^{\rm{c}}$. $K$ and $\{\xi:\, |\xi|=1\}$ can be covered by a finite number of neighborhoods $\{U(x_i)\}_{i=1}^N$ and $\{\Gamma_{x_i}(\xi_j)\}_{j=1}^{M(y_i)}$ respectively, such that on each $U(x_i)\times\Gamma_{x_i}(\xi_j)$ the estimate 
\[
|\nabla_\xi\phi_\eps(y,\xi)|^2\ge c_{i,j}\eps^{r_{i,j}}
\]
holds for some constants $c_{i,j}>0, r_{i,j}\in\R$ and for all $\eps\in(0,\eta_{i,j}]$. It follows that there exist $c>0$, $r\in\R$ and $\eta\in(0,1]$ such that when $y$ is varying in $\cup_{i=1}^N U(x_i)$, $\xi\in\R^p\setminus 0$ and $\eps\in(0,\eta]$ we have 
\[
|\nabla_\xi\phi_\eps(y,\xi)|^2\ge c\eps^r.
\]
This proves that $|\nabla_\xi\phi|^2$ is invertible.
\end{proof}
\begin{remark}
\label{rem_phase_func}
Proposition \ref{prop_R_phi} says that on every relatively compact open subset $U$ of $R_\phi$, $\phi|_{U\times\R^p}$ has the property of a generalized phase function in $\wt{\Phi}(U;\R^p)$. More precisely $\phi_\eps|_{U\times\R^p}\in\Phi(U;\R^p)$ when $\eps$ is varying in some smaller interval $(0,\eta]\subseteq(0,1]$, the net $(\phi_\eps|_{U\times\R^p})_{\eps\in(0,\eta]}$ satifies the ${S^1_{\rm{hg}}(U\times\R^p\setminus 0)}$-moderateness condition and $(\nabla_\xi\phi_\eps(y,\xi/|\xi|))_\eps$ is strictly nonzero. In order to have a representative of $\phi$ defined on the interval $(0,1]$ we may take 
\[
\begin{cases} {\phi_\eps}(y,\xi) & (y,\xi)\in\Om_0\times\R^p,\ \eps\in(0,\eta]\\
\xi & (y,\xi)\in\Om_0\times\R^p,\ \eps\in(\eta,1].
\end{cases}
\]
Clearly the phase function $\phi|_{U\times\R^p}$ does not depend on the choice of the classical phase function we use on the interval $(\eta,1]$.
\end{remark}
The more specific assumption of slow scale-invertibility concerning the generalized symbol $|\nabla_\xi\phi|^2$ is employed in the definition of the following sets.
\begin{definition}
\label{def_C_phi_ssc}
Let $\phi\in\wt{\Phi}(\Om\times\R^p)$. We define $C^\ssc_{\phi}\subseteq\Om\times\R^p\setminus 0$ as the complement of the set of all $(x_0,\xi_0)\in\Om\times\R^p\setminus 0$ with the property that there exist a relatively compact open neighborhood $U(x_0)$ of $x_0$ and a conic open neighborhood $\Gamma(\xi_0)\subseteq\R^p\setminus 0$ of $\xi_0$ such that $|\nabla_\xi\phi|^2$ is slow scale-invertible on $U(x_0)\times\Gamma(\xi_0)$. We set $\pi_\Om(C^\ssc_{\phi})=S^\ssc_{\phi}$ and $R^\ssc_{\phi}=(S^\ssc_{\phi})^{{\rm{c}}}$. 
\end{definition}
By construction $C^\ssc_{\phi}$ is a conic closed subset of $\Om\times\R^p\setminus 0$ and $R^\ssc_{\phi}\subseteq R_{\phi}\subseteq\Om$ is open. In analogy with Proposition \ref{prop_R_phi} we can prove that $|\nabla_\xi\phi|^2$ is slow scale-invertible on $R^\ssc_\phi\times\R^p\setminus 0$.  
\begin{theorem}
\label{theorem_R_phi}
Let $\phi\in\wt{\Phi}(\Om\times\R^p)$ and $a\in\wt{\mathcal{S}}^m_{\rho,\delta}(\Om\times\R^p)$. 
\begin{itemize}
\item[(i)] The restriction $I_{\phi}(a)|_{R_\phi}$ of the functional $I_\phi(a)$ to the region $R_\phi$ belongs to $\G(R_\phi)$.
\item[(ii)] If $\phi\in\wt{\Phi}_\ssc(\Om\times\R^p)$ and $a\in\wt{\mathcal{S}}^m_{\rho,\delta,\rm{rg}}(\Om\times\R^p)$ then $I_{\phi}(a)|_{{R^\ssc_\phi}}\in\Ginf(R^\ssc_\phi)$.
\end{itemize}
\end{theorem}
\begin{proof}
$(i)$ Let $\Om_0$ be a relatively compact open subset of $R_\phi$. By Remark \ref{rem_phase_func} we know that $\phi_0:=\phi|_{\Om_0\times\R^p}$ is a generalized phase function in $\wt{\Phi}(\Om_0;\R^p)$. By Theorem \ref{theorem_map} we have that the oscillatory integral  
\[
\int_{\R^p}e^{i\phi_{0}(y,\xi)}{a}(y,\xi)\, \dslash\xi
\]
defines a generalized function $w_0$ in $\G(\Om_0)$. 
Let now $(\Om_j)_{j\in\N}$ be an open covering of $R_\phi$ such that each $\Om_j$ is relatively compact.
Arguing as above we obtain a sequence of generalized phase functions $\phi_j\in\wt{\Phi}(\Om_j;\R^p)$ and a coherent sequence of generalized functions 
\[
w_j(y)=\int_{\R^p}e^{i\phi_j(y,\xi)}a(y,\xi)\, \dslash\xi\in\G(\Om_j).
\]
Thus the sheaf property of $\G(R_\phi)$ yields the existence of a unique $w\in\G({R_\phi})$ such that $w|_{\Om_j}=w_j$ for all $j$. It remains to prove that
\[
I_{\phi}(a)(u)=\int_{R_\phi}w(y)u(y)\, dy
\]
for $u\in\Gc(R_\phi)$.
On the level of representatives, we may assume that $\supp\, u_\eps$ is contained in some $\Omega_j$ for all $\eps$ and we write the oscillatory integral $\int_{\Om_j\times\R^p}e^{i\phi_{j,\eps}(y,\xi)}a_\eps(y,\xi)u_\eps(y) dy\dslash\xi$ as an iterated one. This yields the following equality between equivalence classes:
\[
I_\phi(a)(u)=\int_{\Om_j\times\R^p}e^{i\phi_j(y,\xi)}a(y,\xi)u(y)\, dy\, \dslash\xi =\int_{\Om_j}w_j(y)u(y)\, dy =\int_{R_\phi}w(y)u(y)\, dy.
\]
$(ii)$ Let $\phi\in\wt{\Phi}_\ssc(\Om\times\R^p)$ and $a\in\wt{\mathcal{S}}_{\rho,\delta,\rm{rg}}^m(\Om\times\R^p)$. Then one easily sees that $\phi_j\in\wt{\Phi}_{\ssc}(\Om_j;\R^p)$ for all $j$ and that by Theorem \ref{theorem_Ginf_map} each $w_j$ is a regular generalized function. Hence $w\in \Ginf(R^\ssc_\phi)$ or, in other words, $I_{\phi}(a)|_{{R^\ssc_\phi}}$ belongs to $\Ginf(R^\ssc_\phi)$.
\end{proof}
Theorem \ref{theorem_R_phi} means that
\[
\singsupp_\G\, I_\phi(a)\subseteq S_\phi
\]
if $\phi\in\wt{\Phi}(\Om\times\R^p)$ and $a\in\wt{\mathcal{S}}^m_{\rho,\delta}(\Om\times\R^p)$ and that
\[
\singsupp_{\Ginf} I_\phi(a)\subseteq S^\ssc_\phi
\]
if $\phi\in\wt{\Phi}_\ssc(\Om\times\R^p)$ and $a\in\wt{\mathcal{S}}^m_{\rho,\delta,\rm{rg}}(\Om\times\R^p)$. 
\begin{example}\label{G_example2}
Returning to Example \ref{G_example1} we are now in the position to analyze the regularity
properties of the generalized kernel functional $I_\phi(a)$ of the solution operator
$A$ corresponding to the hyperbolic Cauchy-problem. For any $v\in \Gc(\R^3)$ we have
\begin{equation}\label{hypsol}
  I_\phi(a) (v) = \iiint e^{i\phi(x,t,y,\xi)}\; a(x,t,y,\xi)\, v(x,t,y)\, dx\, dt\, dy\, \dslash\xi,
\end{equation}
where $a$ and $\phi$ are as in Example \ref{G_example1}. Note that in the case of partial differential operators with smooth coefficients and
distributional initial values the wave front set of the distributional kernel of $A$
determines the propagation of singularities from the initial data. When the
coefficients are non-differentiable functions, or even distributions or generalized
functions, matters are not yet understood in sufficient generality. Nevertheless, the
above results allow us to identify regions where the generalized kernel functional
agrees with a generalized function or is even guaranteed to be a $\Ginf$-regular
generalized function. To identify the set $C_\phi$ in the situation of Example
\ref{G_example1} one simply has to study invertibility of $\d_\xi \phi (x,t,y,\xi) =
\gamma(x,t;0)- y$ as a generalized function in a neighborhood of any given point
$(x_0,t_0,y_0)$.

Under the assumptions on $c$ of Example \ref{G_example1}, the representing nets
$(\gamma_\eps(.,.;0))_{\eps\in(0,1]}$ of $\gamma$ are uniformly bounded on compact sets (e.g., when $c$ is
a bounded generalized constant). For given $(x_0,t_0)$ define the generalized domain of
dependence $D(x_0,t_0)\subseteq\R$ to be the set of accumulation points of the net
$(\gamma_\eps(x_0,t_0;0))_{\eps\in(0,1]}$. Then we have that
$$
   \{(x_0,t_0,y_0) \in\R^3 : y_0 \not\in D(x_0,t_0) \} \subseteq R_\phi.
$$

When $c\in\wt{\R}$ this may be proved by showing that if $(x_0,t_0,y_0)\in C_\phi$ then there exists an accumulation point $c'$ of a representative $(c_\eps)_\eps$ of $c$ such that $y_0=x_0-c't_0$. 
\end{example}

\begin{example}
\label{example_phase_reg}
As an illustrative example concerning the regions involving the regularity of the functional $I_\phi(a)$ we consider the generalized phase function on $\R^2\times\R^2$ given by $\phi_\eps(y_1,y_2,\xi_1,\xi_2)=-\eps y_1\xi_1-s_\eps y_2\xi_2$ where $(s_\eps)_\eps$ is bounded and $(s_\eps^{-1})_\eps$ is a slow scale net. Clearly $\phi:=[(\phi_\eps)_\eps]\in\wt{\Phi}_\ssc(\R^2\times\R^2)$. Simple computations show that $R_\phi=\R^2\setminus(0,0)$ and $R^\ssc_\phi=\R^2\setminus\{y_2=0\}$. We leave it to the reader to check that the oscillatory integral
\[
\int_{\R^2}e^{i\phi(y,\xi)}(1+\xi_1^2+\xi_2^2)^{\frac{1}{2}}\, \dslash\xi =\biggl[\biggl(\int_{\R^2}e^{-i\eps y_1\xi_1-is_\eps y_2\xi_2}(1+\xi_1^2+\xi_2^2)^{\frac{1}{2}}\, \dslash\xi_1\, \dslash\xi_2\biggr)_\eps\biggr]
\]
defines a generalized function in $\R^2\setminus(0,0)$ whose restriction to $\R^2\setminus\{y_2=0\}$ is regular.
\end{example}

The Colombeau-regularity of the functional $I_\phi(a)$ is easily proved in the case of generalized symbols of order $-\infty$.  
\begin{proposition}
\label{prop_smooth_sing}
\begin{itemize}
\item[{\, }]
\item[(i)] If $\phi\in\wt{\Phi}(\Om\times\R^p)$ and $a\in\wt{\mathcal{S}}^{-\infty}(\Om\times\R^p)$ then $\singsupp_{\G}I_\phi(a)=\emptyset$.
\item[(ii)] If $\phi\in\wt{\Phi}_\ssc(\Om\times\R^p)$ and $a\in\wt{\mathcal{S}}^{-\infty}_{\rm{rg}}(\Om\times\R^p)$ then $\singsupp_{\Ginf}I_\phi(a)=\emptyset$.
\end{itemize}
\end{proposition}
\begin{proof}
$(i)$ Arguing as in Section 2 we have that $e^{i\phi(y,\xi)}a(y,\xi)$ is a well-defined element of $\wt{\mathcal{S}}^{-\infty}(\Om\times\R^p)$. Hence by Proposition \ref{pro_C_k}$(ii)$ we obtain that $\int_{\R^p}e^{i\phi(y,\xi)}a(y,\xi)\, \dslash\xi\in\G(\Om)$. A direct inspection of the action of $I_\phi(a)$ at the level of representatives shows that the functional $I_\phi(a)$ coincides with the generalized function $\int_{\R^p}e^{i\phi(y,\xi)}a(y,\xi)\, \dslash\xi$. This means that $\singsupp_\G\, I_\phi(a)=\emptyset$.

$(ii)$ When $\phi$ is a slow scale generalized phase function and $a\in\wt{\mathcal{S}}^{-\infty}_{\rm{rg}}(\Om\times\R^p)$ then $e^{i\phi(y,\xi)}a(y,\xi)\in\wt{\mathcal{S}}^{-\infty}_{\rm{rg}}(\Om\times\R^p)$. Therefore, by Proposition \ref{pro_C_k}$(iii)$ we have that $\int_{\R^p}e^{i\phi(y,\xi)}a(y,\xi)\, \dslash\xi\in\Ginf(\Om)$ and then $\singsupp_{\Ginf} I_\phi(a)=\emptyset$.
\end{proof}

For technical reasons, we will multiply the generalized symbol $a\in\wt{\mathcal{S}}^m_{\rho,\delta}(\Om\times\R^p)$ by a cut-off function $p\in\Cinf(\R^p)$ such that $p(\xi)=0$ for $|\xi|\le 1$ and $p(\xi)=1$ for $|\xi|\ge 2$ in the sequel. 
One easily sees that $a(y,\xi)p(\xi)\in\wt{\mathcal{S}}^m_{\rho,\delta}(\Om\times\R^p)$. Let now $V$ be a closed conic neighborhood of ${\rm{cone\, supp}}\, a$. There exists a smooth function $\chi(y,\xi)\in\Om\times\R^p\setminus 0$, homogeneous of degree 0 in $\xi$ such that $\supp\, \chi\subseteq V$ and $\chi$ is identically $1$ in a smaller neighborhood of ${\rm{cone\, supp}}\, a$ when $|\xi|>1$. It follows that $p(\xi)a(y,\xi)\chi(y,\xi)=a(y,\xi)p(\xi)$ in $\wt{\mathcal{S}}^m_{\rho,\delta}(\Om\times\R^p)$. Note that the representing net $(p(\xi)a_\eps(y,\xi)\chi(y,\xi))_\eps$ of $p(\xi)a(y,\xi)$ is supported in the conic neighborhood $V$ of ${\rm{cone\, supp}}\, a$ uniformly with respect to $\eps\in(0,1]$. 

Before stating Proposition \ref{prop_smooth_sing_cone} we note that if $x_0\in(\pi_\Om(C_\phi\cap {\rm{cone\,supp}}\, a))^{\rm{c}}$ then there exists a relatively compact open neighborhood $U(x_0)$ of $x_0$ and a closed conic neighborhood $V$ of ${\rm{cone\, supp}}\, a$ such that 
\[
(\overline{U(x_0)}\times\{\xi:\, |\xi|\ge 1\})\cap C_\phi\cap V=\emptyset .
\]
\begin{proposition}
\label{prop_smooth_sing_cone}
\begin{itemize}
\item[{\, }]
\item[(i)] If $\phi\in\wt{\Phi}(\Om\times\R^p)$ and $a\in\wt{\mathcal{S}}^m_{\rho,\delta}(\Om\times\R^p)$ then $$\singsupp_{\G}\,I_\phi(a)\subseteq\pi_\Om(C_\phi\cap {\rm{cone\, supp}}\, a).$$
\item[(ii)] If $\phi\in\wt{\Phi}_\ssc(\Om\times\R^p)$ and $a\in\wt{\mathcal{S}}^{m}_{\rho,\delta,\rm{rg}}(\Om\times\R^p)$ then $$\singsupp_{\Ginf}I_\phi(a)\subseteq\pi_\Om(C^\ssc_\phi\cap {\rm{cone\, supp}}\, a).$$
\end{itemize}
\end{proposition}
\begin{proof}
$(i)$ Since $(1-p(\xi))a(y,\xi)\in\wt{\mathcal{S}}^{-\infty}(\Om\times\R^p)$, Proposition \ref{prop_smooth_sing}$(i)$ yields that  $\singsupp_{\G}\,I_\phi(a)$ coincides with $\singsupp_{\G}\,I_\phi(ap)$. By the previous considerations on the conic support of $a$ we can insert a cut-off function $\chi$ with support contained in a closed conic neighborhood $V$ of ${\rm{cone\, supp}}\, a$ as above. Hence, $\singsupp_{\G}\,I_\phi(a)=\singsupp_{\G}\,I_\phi(ap\chi)$. Let now $x_0$ be a point of $(\pi_\Om(C_\phi\cap {\rm{cone\,supp}}\, a))^{\rm{c}}$ and $U(x_0)$ a relatively compact open neighborhood of $x_0$ such that $(\overline{U(x_0)}\times\{\xi:\, |\xi|\ge 1\})\cap C_\phi\cap V=\emptyset$. The generalized symbol $ap\chi$, when restricted to the region $U(x_0)$, has the representative $(a_\eps(y,\xi)p(\xi)\chi(y,\xi))$ which is identically $0$ on $(U(x_0)\times\R^p)\cap C_\phi$. By Theorem \ref{theorem_R_phi}$(i)$ this means that the oscillatory integral $(\int_{\R^p}e^{i\phi(y,\xi)}a(y,\xi)p(\xi)\chi(y,\xi)\, \dslash\xi)|_{{U(x_0)}}$ defines a generalized function in $\G(U(x_0))$ whose representatives can be written in the form 
\[
\int_{\R^p}e^{i\phi_\eps(y,\xi)}L^k_{\phi_\eps}(y;\xi,\partial_\xi)(a_\eps(y,\xi)p(\xi)\chi(y,\xi))\, \dslash\xi
\]
for $\eps$ small enough. This proves that $x_0\not\in\singsupp_{\G}\,I_\phi(a)$.

$(ii)$ If $a\in\wt{\mathcal{S}}^{m}_{\rho,\delta,\rm{rg}}(\Om\times\R^p)$ then $(1-p(\xi))a(y,\xi)\in\wt{\mathcal{S}}^{-\infty}_{\rm{rg}}(\Om\times\R^p)$. By Proposition \ref{prop_smooth_sing}$(ii)$ it follows that $\singsupp_{\Ginf} I_\phi(a)=\singsupp_{\Ginf} I_\phi(ap\chi)$. Arguing as in Theorem \ref{theorem_R_phi}$(ii)$ we obtain that if $x_0\in(\pi_\Om(C_\phi\cap {\rm{cone\,supp}}\, a))^{\rm{c}}$ then $(\int_{\R^p}e^{i\phi(y,\xi)}a(y,\xi)p(\xi)\chi(y,\xi)\, \dslash\xi)|_{{U(x_0)}}$ belongs to $\Ginf(U(x_0))$. Consequently $x_0\not\in\singsupp_{\Ginf}\,I_\phi(a)$.
\end{proof}
\begin{remark}
\label{rem_ess}
When we deal with generalized symbols of refined order the conic support can be substituted by the microsupport. More precisely a combination of Proposition \ref{prop_smooth_sing} with Proposition \ref{prop_smooth_sing_cone} yields the following assertions:
\begin{itemize}
\item[(i)] if $\phi\in\wt{\Phi}(\Om\times\R^n)$ and $a\in\wt{\mathcal{S}}^{\,m/-\infty}_{\rho,\delta}(\Om\times\R^n)$ then $$\singsupp_{\G}\,I_\phi(a)\subseteq\pi_\Om(C_\phi\cap \mu_\G(a));$$
\item[(ii)] if $\phi\in\wt{\Phi}_\ssc(\Om\times\R^n)$ and $a\in\wt{\mathcal{S}}^{\,m/-\infty}_{\rho,\delta,\rm{rg}}(\Om\times\R^n)$ then $$\singsupp_{\Ginf}I_\phi(a)\subseteq\pi_\Om(C^\ssc_\phi\cap \mu_{\Ginf}(a)).$$
\end{itemize}
Indeed, assuming that the cut-off $\chi$ is identically $1$ in a conic neighborhood of $\mu_\G(a)$ when $|\xi|\ge 1$ then $ap(1-\chi)\in\wt{\mathcal{S}}^{-\infty}(\Om\times\R^n)$. Hence $\singsupp_\G I_\phi(a)=\singsupp_\G I_\phi(ap\chi)\subseteq\pi_\Om(C_\phi\cap {\rm{cone\, supp}}(ap\chi))\subseteq\pi_\Om(C_\phi\cap \mu_\G(a))$. By means of analogous arguments one can easily prove the second inclusion above for $\phi\in\wt{\Phi}_\ssc(\Om\times\R^n)$ and $a\in\wt{\mathcal{S}}^{\,m/-\infty}_{\rho,\delta,\rm{rg}}(\Om\times\R^n)$.
\end{remark}

We conclude the paper by investigating the $\G$-wave front set and the $\Ginf$-wave front set of the functional $I_\phi(a)$ under suitable assumptions on the generalized symbol $a$ and the phase function $\phi$. We leave it to the reader to check that when $\phi\in\wt{\Phi}(\Om\times\R^p)$, $U\subseteq\overline{U}\Subset\Om$, $\Gamma\subseteq\R^n\setminus 0$, $V\subseteq \Om\times\R^p\setminus 0$, then 
\[
\mathop{{\rm{Inf}}}\limits_{\substack{y\in U, \xi\in \Gamma\\ (y,\theta)\in V}}\frac{|\xi-\nabla_y\phi(y,\theta)|}{|\xi|+|\theta|}:=\biggl[\biggl(\inf_{\substack{y\in U, \xi\in \Gamma\\ (y,\theta)\in V}}\frac{|\xi-\nabla_y\phi_\eps(y,\theta)|}{|\xi|+|\theta|}\biggr)_\eps\biggr]
\]
is a well-defined element of $\wt{\C}$. 
\begin{theorem}
\label{theorem_WF}
\leavevmode
\begin{itemize}
\item[(i)] Let $\phi\in\wt{\Phi}(\Om\times\R^p)$ and $a\in\wt{\mathcal{S}}^m_{\rho,\delta}(\Om\times\R^p)$. The generalized wave front set $\WF_\G I_\phi(a)$ is contained in the set $W_{\phi,a}$ of all points $(x_0,\xi_0)\in\CO{\Om}$ with the property that for all relatively compact open neighborhoods $U(x_0)$ of $x_0$, for all open conic neighborhoods $\Gamma(\xi_0)\subseteq \R^n\setminus 0$ of $\xi_0$, for all open conic neighborhoods $V$ of {\rm{cone\,supp}}\,$a\cap C_\phi$ such that $V\cap (U(x_0)\times\R^p\setminus 0)\neq\emptyset$ the generalized number 
\beq
\label{gen_num_inv}
\mathop{{\rm{Inf}}}\limits_{\substack{y\in U(x_0), \xi\in \Gamma(\xi_0)\\ (y,\theta)\in V\cap( U(x_0)\times\R^p\setminus 0)}}\frac{|\xi-\nabla_y\phi(y,\theta)|}{|\xi|+|\theta|}
\eeq
is not invertible.
\item[(ii)] If $\phi\in\wt{\Phi}_\ssc(\Om\times\R^p)$ and $a\in\wt{\mathcal{S}}^m_{\rho,\delta,\rm{rg}}(\Om\times\R^p)$ then $\WF_{\Ginf}I_\phi(a)$ is contained in the set $W^\ssc_{\phi,a}$ of all points $(x_0,\xi_0)\in\CO{\Om}$ with the property that for all relatively compact open neighborhoods $U(x_0)$ of $x_0$, for all open conic neighborhoods $\Gamma(\xi_0)\subseteq \R^n\setminus 0$ of $\xi_0$, for all open conic neighborhoods $V$ of {\rm{cone\,supp}}\,$a\cap C^\ssc_\phi$ such that $V\cap (U(x_0)\times\R^p\setminus 0)\neq\emptyset$ the generalized number \eqref {gen_num_inv}
is not slow scale-invertible.
\end{itemize}
\end{theorem}
\begin{proof}
$(i)$ Let us assume that $(x_0,\xi_0)\not\in W_{\phi,a}$. By the definition of $W_{\phi,a}$ there exist an open relatively compact neighborhood $U(x_0)$ of $x_0$, an open conic neighborhood $\Gamma(\xi_0)$ of $\xi_0$ and a closed conic neighborhood $V$ of ${\rm{cone\,supp}}\,a\cap C_\phi$ such that the corresponding generalized number in \eqref{gen_num_inv} is invertible. Assume that $\chi(y,\theta)$ is a smooth function on $\Om\times\R^p\setminus 0$ , homogeneous of degree $0$ in $\theta$ with support contained in $V$ and identically $1$ in a smaller neighborhood $V'$ of ${\rm{cone\, supp}}\,a\cap C_\phi$ when $|\theta|\ge 1$. Choosing $p(\theta)$ as explained before Proposition \ref{prop_smooth_sing_cone}, we have that $I_\phi(a(1-p))\in\G(\Om)$ by Proposition \ref{prop_smooth_sing}$(i)$. Hence $\WF_\G I_\phi(a)= \WF_\G I_\phi(ap)$. Inserting the cut-off function $\chi$ and making use of Proposition \ref{prop_smooth_sing_cone}$(i)$ we arrive at $$\singsupp_\G\, I_\phi(ap(1-\chi))\subseteq \pi_\Om (C_\phi\cap {\rm{cone\, supp}}a\cap (V')^{\rm{c}})=\emptyset.$$ Thus, $\WF_\G I_\phi(a)=\WF_\G I_\phi(ap\chi)$. 

Since $(x_0,\xi_0)\not\in W_{\phi,a}$ there exists a representative $(\phi_\eps)_\eps$ of $\phi$, $\eta\in(0,1]$ and $r\in\R$ such that 
\beq
\label{est_below}
|\xi-\nabla_y\phi_\eps(y,\theta)|\ge \eps^r(|\xi|+|\theta|)
\eeq
for all $y\in U(x_0)$, $\xi\in\Gamma(\xi_0)$ and $(y,\theta)\in V\cap U(x_0)\times\R^p\setminus 0$. Consider the family of differential operators
\[
L_\eps:=\sum_{j=1}^n\frac{\xi_j-\partial_{y_j}\phi_\eps(y,\theta)}{|\xi-\nabla_y\phi_\eps(y,\theta)|^2}D_{y_j}
\]
under the assumptions on $y$, $\theta$, $\xi$ above and denote the coefficient $(\xi_j-\partial_{y_j}\phi_\eps(y,\theta))/|\xi-\nabla_y\phi_\eps(y,\theta)|^2$ by $d_{j,\eps}(y,\theta,\xi)$. Combining \eqref{est_below} with the fact that $(\phi_\eps)_\eps\in\mM_{S^1_{\rm{hg}}(\Om\times\R^p\setminus 0)}$ we get that 
\begin{multline}
\label{est_up}
\forall\alpha\in\N^n\, \exists N\in\N\, \exists\eta\in(0,1]\, \forall \xi\in\Gamma(\xi_0)\, \forall(y,\theta)\in V\cap(U(x_0)\times\R^p\setminus 0)\, \forall\eps\in(0,\eta]\\
|\partial^\alpha_y d_{j,\eps}(y,\theta,\xi)|\le \eps^{-N}(|\theta|+|\xi|)^{-1}
\end{multline}
By construction, $L_\eps e^{i(\phi_\eps(y,\theta)-y\xi)}=e^{i(\phi_\eps(y,\theta)-y\xi)}$ and the transpose operator is of the form
\[
{\ }^tL_\eps=\sum_{j=1}^n a_{j,\eps}(y,\theta,\xi)D_{y_j}+c_\eps(y,\theta,\xi),
\]
with the coefficients $(a_{j,\eps})_\eps$ and $(c_\eps)_\eps$ satisfying condition \eqref{est_up}. If  $(b_\eps)_\eps\in\mM_{S^m_{\rho,\delta}(\Om\times\R^p)}$, an induction argument yields the following result:
\begin{multline}
\label{est_transp}
\forall l\in\N\, \exists N\in\N\, \exists\eta\in(0,1]\, \forall \xi\in\Gamma(\xi_0)\, \forall(y,\theta)\in  V\cap(U(x_0)\times\R^p\setminus 0)\, \forall\eps\in(0,\eta]\\
|({\ }^tL_\eps)^l b_\eps(y,\theta)|\le \eps^{-N}(|\xi|+|\theta|)^{-l}(1+|\theta|)^{m+\delta l}.
\end{multline} 
We now have all the tools at hand for dealing with the Fourier transform of the functional ${\varphi I_\phi(ap\chi)}\in\LL(\G(\Om),\wt{\C})$ when $\varphi\in\Cinfc(U)$. By definition, it is the tempered generalized function given by the integral 
\[
\int_{\Om\times\R^p}e^{i\phi(y,\theta)}a(y,\theta)p(\theta)\chi(y,\theta)e^{-iy\xi}\varphi(y)\, dy\, \dslash\theta,
\]
and from the previous considerations it follows that it has a representative $(f_\eps)_\eps$ of the form 
\[
\int_{\Om\times\R^p}e^{i(\phi_\eps(y,\theta)-y\xi)}({\ }^tL_\eps)^l (a_\eps(y,\theta)p(\theta)\chi(y,\theta)\varphi(y))\, dy\, \dslash\theta
\]
when $\eps$ is small enough and $\xi$ is varying in $\Gamma(\xi_0)$. Making use of the estimate \eqref{est_transp} and taking $l$ so large that for $\delta<\lambda<1$ one has $m+(\delta -\lambda)l<-p$, we conclude that 
\begin{multline}
\label{est_fourier}
\exists N_l\in\N\, \exists c>0\, \exists\eta\in(0,1]\, \forall\xi\in\Gamma(\xi_0)\cap\{\xi: |\xi|\ge 1\}\, \forall\eps\in(0,\eta]\\
|f_\eps(\xi)|\le c\eps^{-N_l}|\xi|^{-l(1-\lambda)}\int_{\R^p}(1+|\theta|)^{m+(\delta-\lambda)l}\, \dslash\theta.
\end{multline}
This shows that $(\varphi I_\phi(ap\chi))^{\widehat{\, }}$ is a generalized function in $\G_{\S\hskip-2pt,0}(\Gamma)$. The characterization of the wave front set of a functional given in \cite{Garetto:06a} proves that $(x_0,\xi_0)$ does not belong to $\WF_\G I_\phi(ap\chi)=\WF_\G I_\phi(a)$.

$(ii)$ We now work with $\phi\in\wt{\Phi}_\ssc(\Om\times\R^p)$, $a\in\wt{\mathcal{S}}^{m}_{\rho,\delta,\rm{rg}}(\Om\times\R^p)$ and $(x_0,\xi_0)\not\in W^\ssc_{\phi,a}$. Choosing $p(\theta)$ and $\chi(x,\theta)$ as in the first case with $V$ and $V'$ neighborhoods of ${\rm{cone\, supp}}\,a\cap C^\ssc_\phi$, we have that $I_\phi(a(1-p))\in\Ginf(\Om)$ by Proposition \ref{prop_smooth_sing}$(ii)$. Moreover, Proposition \ref{prop_smooth_sing_cone}$(ii)$ leads to $$\singsupp_\G\, I_\phi(ap(1-\chi))\subseteq \pi_\Om (C^\ssc_\phi\cap {\rm{cone\, supp}}a\cap (V')^{\rm{c}})=\emptyset.$$ Hence, $\WF_{\Ginf} I_\phi(a)=\WF_{\Ginf} I_\phi(ap\chi)$. By the definition of $W_{\phi,a}^\ssc$ there exists a representative $(\phi_\eps)_\eps$ of $\phi$, a slow scale net $(s_\eps)_\eps$ and a number $\eta\in(0,1]$ such that 
\beq
\label{est_below_ssc}
|\xi-\nabla_y\phi_\eps(y,\theta)|\ge s_\eps^{-1}(|\xi|+|\theta|)
\eeq
for all $\xi\in\Gamma(\xi_0)$, $(y,\theta)\in V\cap (U(x_0)\times\R^p\setminus 0)$ and $\eps\in(0,\eta]$. This combined with the hypothesis $\phi\in\wt{\Phi}_\ssc(\Om\times\R^p)$ implies that the coefficients of the operator ${\ }^tL_\eps$ are nets of slow scale type. Further, when $(b_\eps)_\eps$ is the representative of a generalized symbol $b$ in $\wt{\mathcal{S}}^m_{\rho,\delta,\rm{rg}}(\Om\times\R^p)$, we are allowed to change the order of the quantifiers $\forall l\in\N$ and $\exists N\in\N$ in \eqref{est_transp}. As a consequence, the estimate \eqref{est_fourier} holds with some $N$ independent of $l\ge 0$. We conclude that $(\varphi I_\phi(ap\chi))^{\widehat{\, }}$ is a generalized function in $\Ginf_{\S\hskip-2pt,0}(\Gamma)$ and then $(x_0,\xi_0)\not\in\WF_{\Ginf} I_\phi(ap\chi)=\WF_{\Ginf} I_\phi(a)$.
\end{proof}
\begin{example}
\label{C_example3}
Theorem \ref{theorem_WF} can be employed for investigating the generalized wave front sets of the kernel $K_A:=I_\phi(a)$ of the Fourier integral operator introduced in Example \ref{G_example1}. For simplicity we assume that $c$ is a bounded generalized constant in $\wt{\R}$ and that $a=1$. Let $((x_0,t_0,y_0),\xi_0)\in\WF_\G K_A$. From the first assertion of Theorem \ref{theorem_WF} we know that the generalized number given by 
\beq
\label{non_inv_ex}
\inf_{\substack{(x,t,y)\in U, \xi\in \Gamma\\ ((x,t,y),\theta)\in V\cap(U\times\R\setminus 0)}}\frac{|\xi-(\theta,-c_\eps\theta,-\theta)|}{|\xi|+|\theta|}
\eeq
is not invertible, for every choice of neighborhoods $U$ of $(x_0,t_0,y_0)$, $\Gamma$ of $\xi_0$ and $V$ of $C_\phi$. Note that it is not restrictive to assume that $|\theta|=1$. We fix some sequences $(U_n)_n$, $(\Gamma_n)_n$ and $(V_n)_n$ of neighborhoods shrinking to $(x_0,t_0,y_0)$, $\{\xi_0\lambda:\lambda>0\}$ and $C_\phi$ respectively. By \eqref{non_inv_ex} we find a sequence $\eps_n$ tending to $0$ such that for all $n\in\N$ there exists $\xi_n\in\Gamma_n$, $(x_n,t_n,y_n,\theta_n)\in V_n$ with $|\theta_n|=1$ and $(x_n,t_n,y_n)\in U_n$ such that
\[
|\xi_n-(\theta_n,-c_{\eps_n}\theta_n,-\theta_n)|\le \eps_n(|\xi_n|+1).
\]
In particular, $\xi_n$ remains bounded. Passing to suitable subsequences we obtain that there exist $\theta$ such that $(x_0,t_0,y_0,\theta)\in C_\phi$, an accumulation point $c'$ of $(c_\eps)_\eps$ and a multiple $\xi'$ of $\xi_0$ such that $\xi'=(\theta,-c'\theta,-\theta)$. It follows that
\[
\frac{\xi_0}{|\xi_0|}=
\frac{\xi'}{|\xi'|}=\frac{1}{\sqrt{2+(c')^2}|\theta|}\,(\theta,-c'\theta,-\theta).
\]

In other words the $\G$-wave front set of the kernel $K_A$ is contained in the set of points of the form $((x_0,t_0,y_0),(\theta_0,-c'\theta_0,-\theta_0))$ where $(x_0,t_0,y_0,\theta_0)\in C_\phi$ and $c'$ is an accumulation point of a net representing $c$. Since in the classical case (when $c\in\R$) the distributional wave front set of the corresponding kernel is the set $\{((x_0,t_0,y_0),(\theta_0,-c\theta_0,-\theta_0)):\, (x_0,t_0,y_0,\theta_0)\in C_\phi\}$, the result obtained above for $\WF_\G K_A$ is a generalization in line with what we deduced about the regions $R_\phi$ and $C_\phi$ in Example \ref{G_example2}.

\end{example}

\begin{remark}
\label{rem_class_phase}
It is instructive to give an interpretation of the results stated in Theorem \ref{theorem_WF} in the case of classical phase functions.
 
When $\phi$ is a classical phase function the set $W_{\phi,a}$ as well as the set $W^\ssc_{\phi,a}$ coincide with 
\beq
\label{set_class}
\{(x,\nabla_x\phi(x,\theta)):\, (x,\theta)\in{\rm{cone\, supp}}\,a\cap C_\phi\}.
\eeq
Indeed, if $(x_0,\xi_0)$ is in the complement of the region defined in \eqref{set_class} then there exists a relatively compact neighborhood $U(x_0)$ of $x_0$, a closed conic neighborhoods $\Gamma(\xi_0)\subseteq \R^n\setminus 0$ of $\xi_0$ and a conic neighborhood $V$ of {\rm{cone\,supp}}\,$a\cap C^\ssc_\phi$ with $V\cap (U(x_0)\times\R^p\setminus 0)\neq\emptyset$ such that $\nabla_x\phi(x,\theta)\not\in \Gamma(\xi_0)$ for all $(x,\theta)\in V$ with $x\in U(x_0)$. By continuity and homogeneity of $\phi$ we conclude that there exists $c>0$ such that 
\[
\inf_{\substack{y\in U, \xi\in \Gamma\\ (y,\theta)\in V}}\frac{|\xi-\nabla_y\phi(y,\theta)|}{|\xi|+|\theta|}>c.
\]
It follows that $(x_0,\xi_0)\not\in W_{\phi,a}$. Clearly if $(x_0,\xi_0)\not\in W_{\phi,a}$ then $(x_0,\xi_0)$ does not belong to the set in \eqref{set_class}.
\end{remark}

When the functional $I_\phi(a)$ is given by a generalized symbol of refined order on $\Om\times\R^n$ a more precise evaluation of the wave front sets $\WF_\G I_\phi(a)$ and $\WF_{\Ginf} I_\phi(a)$ can be obtained by making use of the generalized microsupports $\mu_\G$ and $\mu_{\Ginf}$ instead of the conic support of $a$ in the definition of the sets $W_{\phi,a}$ and $W^\ssc_{\phi,a}$ respectively.
\begin{corollary}
\label{cor_mu}
\leavevmode
\begin{itemize}
\item[(i)] Let $\phi\in\wt{\Phi}(\Om\times\R^n)$ and $a\in\wt{\mathcal{S}}^{\,m/-\infty}_{\rho,\delta}(\Om\times\R^n)$. The generalized wave front set $\WF_\G I_\phi(a)$ is contained in the set $W_{\phi,a,\G}$ of all points $(x_0,\xi_0)\in\CO{\Om}$ with the property that for all relatively compact open neighborhoods $U(x_0)$ of $x_0$, for all open conic neighborhoods $\Gamma(\xi_0)\subseteq \R^n\setminus 0$ of $\xi_0$, for all open conic neighborhoods $V$ of $\mu_\G\, a\cap C_\phi$ such that $V\cap (U(x_0)\times\R^n\setminus 0)\neq\emptyset$ the generalized number 
\beq
\label{gen_num_inv_2}
\mathop{{\rm{Inf}}}\limits_{\substack{y\in U(x_0), \xi\in \Gamma(\xi_0)\\ (y,\theta)\in V\cap( U(x_0)\times\R^n\setminus 0)}}\frac{|\xi-\nabla_y\phi(y,\theta)|}{|\xi|+|\theta|}
\eeq
is not invertible.
\item[(ii)] If $\phi\in\wt{\Phi}_\ssc(\Om\times\R^n)$ and $a\in\wt{\mathcal{S}}^{\,m/-\infty}_{\rho,\delta,\rm{rg}}(\Om\times\R^n)$ then $\WF_{\Ginf}I_\phi(a)$ is contained in the set $W^\ssc_{\phi,a,{\Ginf}}$ of all points $(x_0,\xi_0)\in\CO{\Om}$ with the property that for all relatively compact open neighborhoods $U(x_0)$ of $x_0$, for all open conic neighborhoods $\Gamma(\xi_0)\subseteq \R^n\setminus 0$ of $\xi_0$, for all open conic neighborhoods $V$ of $\mu_{\Ginf}\, a\cap C^\ssc_\phi$ such that $V\cap (U(x_0)\times\R^n\setminus 0)\neq\emptyset$ the generalized number \eqref {gen_num_inv_2}
is not slow scale-invertible.
\end{itemize}
\end{corollary}
Note that since $\mu_\G\, a\subseteq\mu_{\Ginf}\, a\subseteq{\rm{cone\, supp}}\,a$ for any symbol of refined order, we have that $W_{\phi,a,\G}\subseteq W_{\phi,a}$ and $W_{\phi,a,\Ginf}\subseteq W^{\ssc}_{\phi,a}$ by construction.
\begin{proof}[Proof of Corollary \ref{cor_mu}.]
$(i)$ Let $p\in\Cinf(\R^n)$ such that $p(\theta)=0$ for $|\theta|\le 1$ and $p(\theta)=1$ for $|\theta|\ge 2$ and let $\chi(y,\theta)$ be a smooth function in $\Om\times\R^n\setminus 0$ , homogeneous of degree $0$ in $\theta$ with support contained in a neighborhood $V$ of $\mu_\G(a)\cap C_\phi$ and identically $1$ is a smaller neighborhood $V'$ of $\mu_\G(a)\cap C_\phi$ when $|\theta|\ge 1$. Remark \ref{rem_ess}$(i)$ implies that $\WF_\G I_\phi(a)=\WF_\G I_\phi(ap\chi)$. At this point an application of Theorem \ref{theorem_WF}$(i)$ entails the desired inclusion, since ${\rm{cone\, supp}}(ap\chi)\subseteq\supp\, \chi\subseteq V$.

$(ii)$ The second assertion of the theorem is obtained by a combination of Remark \ref{rem_ess}$(ii)$ with the second statement of Theorem \ref{theorem_WF}.
\end{proof}
\begin{remark}
\label{rem_class_phase_2}
When $\phi$ is a classical phase function on $\Om\times\R^n$, Corollary \ref{cor_mu} and the truncation arguments employed in Remark \ref{rem_class_phase} yield the following inclusions:
\beq
\label{incl_1}
\WF_\G I_\phi(a)\subseteq\{(x,\nabla_x\phi(x,\theta)):\, (x,\theta)\in\mu_\G(a)\cap C_\phi\}, 
\eeq
\beq
\label{incl_2}
\WF_{\Ginf} I_\phi(a)\subseteq\{(x,\nabla_x\phi(x,\theta)):\, (x,\theta)\in\mu_{\Ginf}(a)\cap C_\phi\}, 
\eeq
valid for $a\in\wt{\mathcal{S}}^{\,m/-\infty}_{\rho,\delta}(\Om\times\R^n)$ and $a\in\wt{\mathcal{S}}^{\,m/-\infty}_{\rho,\delta,\rm{rg}}(\Om\times\R^n)$ respectively.
\end{remark}

Finally, we consider a generalized pseudodifferential operator $a(x,D)$ on $\Om$ and its kernel $K_{a(x,D)}\in\LL(\Gc(\Om\times\Om),\wt{\C})$. From \eqref{incl_1}, we have that $\WF_\G(K_{a(x,D)})$ is contained in the normal bundle of the diagonal in $\Om\times\Om$ when $a\in\wt{\mathcal{S}}^{m}_{\rho,\delta}(\Om\times\R^n)$. By \eqref{incl_2}, $\WF_{\Ginf}(K_{a(x,D)})$ is a subset of the normal bundle of the diagonal in $\Om\times\Om$ when $a$ is regular. We define the sets
\[
\WF_\G(a(x,D))=\{(x,\xi)\in\CO{\Om}:\, (x,x,\xi,-\xi)\in\WF_\G(K_{a(x,D)})\}
\]
and
\[
\WF_{\Ginf}(a(x,D))=\{(x,\xi)\in\CO{\Om}:\, (x,x,\xi,-\xi)\in\WF_{\Ginf}(K_{a(x,D)})\}.
\]
\begin{proposition}
\label{prop_pseudo}
Let $a(x,D)$ be a generalized pseudodifferential operator.
\begin{itemize}
\item[(i)] If $a\in\wt{\mathcal{S}}^{m}_{\rho,\delta}(\Om\times\R^n)$ then $\WF_\G(a(x,D))\subseteq\mu\, \supp_\G(a)$.
\item[(ii)]If $a\in\wt{\mathcal{S}}^{m}_{\rho,\delta,{\rm{rg}}}(\Om\times\R^n)$ then $\WF_{\Ginf}(a(x,D))\subseteq\mu\, \supp_{\Ginf}(a)$.
\end{itemize}
\end{proposition}
\begin{proof}
$(i)$ Let $(a_\eps)_\eps$ be a representative of $a$ and $\kappa((a_\eps)_\eps)=(a_\eps)_\eps+\Neg^{-\infty}(\Om\times\R^n)$. From \eqref{incl_1} we have that 
\[
\WF_\G(K_{\kappa((a_\eps)_\eps)(x,D)})\subseteq\{(x,x,\xi,-\xi)\in\CO{\Om\times\Om}:\, (x,\xi)\in\mu_\G(\kappa((a_\eps)_\eps))\}.
\]
The intersection over all representatives of $a$ combined with \eqref{ess_G} yields
\[
\WF_\G(K_{a(x,D)})\subseteq\{(x,x,\xi,-\xi)\in\CO{\Om\times\Om}:\, (x,\xi)\in\mu\, \supp_\G(a)\}.
\]
Hence,
\[
\WF_\G(a(x,D))\subseteq \mu\, \supp_\G(a).
\]
$(ii)$ The second assertion of the proposition is obtained as above from \eqref{incl_2} and the equality \eqref{ess_Ginf}.
\end{proof}
In consistency with the notations introduced in \cite[Definition 3.9]{GH:05} and \cite[Definition 3.11]{Garetto:06a} we can define the $\G$-microsupport of a properly supported generalized pseudodifferential operator as 
\[
\Bigmu\,\supp_\G(A):=\bigcap_{\substack{a\in\wt{\mathcal{S}}^m_{\rho,\delta}(\Om\times\R^n)\\ a(x,D)=A}}\mu\,\supp_\G(a)
\]
and the $\Ginf$-microsupport as 
\[
\Bigmu\,\supp_{\Ginf}(A):=\bigcap_{\substack{a\in\wt{\mathcal{S}}^m_{\rho,\delta,{\rm{rg}}}(\Om\times\R^n)\\ a(x,D)=A}}\mu\,\supp_{\Ginf}(a).
\]
From Proposition \ref{prop_pseudo} it follows that $\WF_\G(A)\subseteq\mu\,\supp_\G(A)$. In particular, if $A$ is given by a symbol in $\wt{\mathcal{S}}^{m}_{\rho,\delta,{\rm{rg}}}(\Om\times\R^n)$ then $\WF_{\Ginf}(A)\subseteq\mu\, \supp_{\Ginf}(A)$.

\bibliographystyle{abbrv}
\newcommand{\SortNoop}[1]{}

\end{document}